\documentstyle[12pt]{article}
\newcommand{\qu}{\quad}
\newcommand{\en}{\enspace}
\newcommand{\bi}{\bigskip}
\newcommand{\me}{\medskip}
\newcommand{\no}{\noindent}                    

\newcommand{\bu}{\bigtriangleup}
\newcommand{\be}{\begin{equation}}    
\newcommand{\ee}{\end{equation}}    
\newcommand{\bea}{\begin{eqnarray}}                                              
\newcommand{\eea}{\end{eqnarray}}

\def\theequation{\thesection.\arabic{equation}}

\begin{document}
\sloppy 
\begin{center} 

{\LARGE \bf Zeta Functions, Determinants and Torsion
 
for Open Manifolds}

\bi   

J\"urgen Eichhorn

Institut f\"ur Mathematik und Informatik

Jahnstra\ss e 15 a, 17487 Greifswald, Germany

e-mail: eichhorn@rz.uni-greifswald.de
\end{center}

\bi

\begin{abstract}

\no 
On an open manifold, the spaces of metrics or connections of bounded
geometry, respectively, split into an uncountable number of components.
We show that for a pair of metrics or connections, belonging to the
same component, relative $\zeta$-functions, determinants, torsion 
for pairs of generalized Dirac operators
are well defined.

\end{abstract}

\section{Introduction}    
                      
Many of the most important invariants which are defined for closed 
manifolds don't make sense for open manifolds. Integrals defining e. g.
characteristic numbers in general diverge. The spectrum of elliptic
self-adjoint operators is not purely discrete etc.. One successful
approach is to restrict to bounded geometry, to fix one metric in this
class and to consider relative invariants. Concerning self-adjoint
differential operators associated to geometry, bounded geometry always
implies that their spectrum contains a half line, i. e. it is very far from
being discrete. Hence $\zeta$-functions don't make sense. But if we fix
a component $comp(g_0)$ in the completed space ${\cal M}^{p,r}(I,B_k)$ of
metrics $g$ satisfying the conditions
 
\me 
\begin{center}
\begin{tabular}{cll}
$(I)$ & {} & $r_{inj} (M,g) = \inf_{x \in M} r_{inj}(x) > 0$ , \\ 
{} \\
$(B_k)$ & {} & $|(\nabla^g)^i R^g| \le C_i, \en 0 \le i \le k$ , 
\end{tabular}
\end{center}
\me                                                 

\no  
then we can consider for instance the pair $\bu_q(g), \bu_q(g_0), g \in
comp(g_0), \bu _q$ the Laplace operator acting on $q$-forms. 
Then, for $k \ge r > n+2, p=1$, we show that 

\be
e^{-t\bu_q(g)} - e^{-t\bu_q(g_0)}  
\ee                     

\no
and                                              

\be
\bu_q(g) e^{-t\bu_q(g)} - \bu _q(g_0)e^{-t\bu_q(g_0)}                           
\ee   
                     
\no
are for $t>0$ of trace class, and their trace norms are uniformly 
bounded on compact $t$--intervalls.

This is a consequence of our extended 
approach to generalized Dirac operators, considering the
completed space $C_E^{p,r}(B_k)$ of Clifford connections. 

\no
Let 
$g' \in comp(g) \subset    
{\cal M}^{1,r} (I,B_k), 
\nabla' \in comp(\nabla) \subset C_E^{1,r}(B_k), 
k \ge r > n+2, 
D=D(g,\nabla), D'=D'(g',\nabla')$ 
be the generalized Dirac operators, then

\be
e^{-tD'^2} - e^{-tD^2} ,  
\ee 
                                            
\be
D'e^{-tD'^2} - De^{-tD^2}   
\ee 
                 
\no
are of trace class and their trace norms are uniformly bounded on compact
$t$--intervalls. Assuming additionally $\inf \sigma_e (D^2)>0$, we 
define relative $\zeta$--functions, determinants and torsion in the case
of the Laplace operator.

The paper is organized as follows. In section 2 we present the necessary
facts concerning Clifford bundles, generalized Dirac operators and 
Sobolev spaces. Section 3 is devoted to  spaces of metrics and
connections. Section 4 contains some general heat kernel estimates which
are needed in section 5. We present in section 5 the first essential
step of our approach, proving that for fixed $g$ and variation of the
Clifford connection $\nabla$ to $\nabla'$ the operators $e^{-tD^2}-e^{-tD'^2}$,
$De^{-tD^2}-D'e^{-tD'^2}$ are of trace class and their trace norm is
uniformly bounded on compact $t$--intervalls $[a_0,a_1]$, $a_0>0$.
Section 6 is devoted to the generalization of 5, admitting variation
of the bundle metric and the Clifford structure too. We apply our
results in sections 7 and 8, establishing certain relative index
theorems and defining $\zeta$--functions, determinants and torsion.
In a forthcoming paper we drop the assumption $\inf \sigma_e(D^2)>0$,
define relative $\eta$--functions and present further applications.

\section{Clifford bundles, generalized Dirac operators and Sobolev spaces}
                                                                         
We recall for completeness very shortly the basic properties of 
generalized Dirac operators on open manifolds. Let $(M^n,g)$ be a
Riemannian manifold, $m \in M, Cl(T_mM,g_m)$ the corresponding Clifford
algebra at $m$. $Cl(T_mM,g_m)$ shall be complexified or not, depending
on the other bundles and structure under consideration. A hermitian
vector bundle $E \rightarrow M$ is called a bundle of Clifford modules
if each fibre $E_m$ is a Clifford module over $Cl(T_mM,g_m)$ with
skew symmetric Clifford multiplication. We assume $E$ to be endowed
with a compatible connection $\nabla^E$, i.e. $\nabla^E$ is metric and 

\[ \nabla^E_X (Y \cdot \Phi) = (\nabla^g_X Y) \cdot \Phi + Y \cdot  
(\nabla^E_X \Phi) , \]

$X,Y \in \Gamma (TM), \Phi \in \Gamma(E)$. Then we call the pair 
$(E,\nabla^E)$ a Clifford bundle. The composition 

\[ \Gamma(E) \stackrel{\nabla}{\longrightarrow} \Gamma (T^* M \otimes E) 
   \stackrel{g}{\longrightarrow} \Gamma ( T M \otimes E) 
   \stackrel{\cdot}{\longrightarrow} \Gamma (E) \]        

\no
shall be called the generalized Dirac operator $D$. We have
$D=D(g,E,\nabla)$. If $X_1, \dots X_n$ is an orthonormal basis in $T_mM$ then

\[ D = \sum^n_{i=1} X_i \cdot \nabla^E_{X_i} . \]

\no
$D$ is of first order elliptic, formally self-adjoint and 

\[ D^2 = \bu^E + {\cal R} , \]

\no
where $\bu^E = (\nabla^E)^* \nabla^E$ and $ {\cal R} \in \Gamma(End(E))$ 
is the 
bundle endomorphism

\[ {\cal R} \Phi = \frac{1}{2} \sum^n_{i,j=1} X_i X_j R^E (X_i,X_j) \Phi . \]

\no
Next we recall some associated functional spaces and their properties
if we assume bounded geometry. These facts are contained in [5], [7], [2].

Let $E \rightarrow M$ be a Clifford bundle, $\nabla = \nabla^E$, $D$ the 
generalized Dirac operator. Then we define for $\Phi \in \Gamma(E),
p \ge 1, r \in {\bf Z}, r \ge 0$,  
\begin{eqnarray} 
   |\Phi|_{W^{p,r}} &:=&  \left( 
   \int \sum^r_{i=0} | \nabla^i \Phi |^p_x dvol_x(g)
   \right)^{\frac{1}{p}} , \nonumber \\
   |\Phi|_{H^{p,r}} &:=&  \left( 
   \int \sum^r_{i=0} | D^i \Phi |^p_x dvol_x(g)
   \right)^{\frac{1}{p}} , \nonumber \\  
   W^p_r(E) &:=& \left\{ 
   \Phi \in \Gamma(E) \big| |\Phi|_{W^{p,r}} < \infty
   \right\} , \nonumber \\
   W^{p,r}(E) &:=& \en \mbox{completion of} 
   \en W^p_r \en \mbox{w. r. t.}
   \en | \en |_{W^{p,r}} , \nonumber \\
   H^p_r(E) &:=& \left\{ 
   \Phi \in \Gamma(E) \big| |\Phi|_{H^{p,r}} < \infty
   \right\} , \nonumber \\
   H^{p,r}(E) &:=& \en \mbox{completion of} 
   \en H^p_r \en \mbox{w. r. t.}
   \en | \en |_{H^{p,r}} . \nonumber 
\end{eqnarray}

\no 
In a great part of our consideration we restrict to $p=2$. In this case
we write $W^{2,r} \equiv W^r, H^{2,r} \equiv H^r $ etc.. If $r<0$ then
we set 
 
\begin{eqnarray} 
W^r(E) &:=& \Big( W^{-r}(E) \Big)^* , \nonumber \\    
H^r(E) &:=& \Big( H^{-r}(E) \Big)^* . \nonumber 
\end{eqnarray}

\no
Assume $(M^n,g)$ complete. Then $C^\infty_c(E)$ is a dense subset of 
$W^{p,1}(E)$ and $H^{p,1}(E)$. This follows from proposition 1.4 in [2].
If we use this density and the fact

\[ |D \Phi(m)| \le C \cdot |\nabla \Phi(m)| , \]

\no 
we obtain $|\Phi|_{H^{p,1}} \le C' \cdot |\Phi|_{W^{p,1}}$ and a
continous embedding

\[ W^{p,1}(E) \hookrightarrow H^{p,1}(E) . \]

\no 
For $r>1$ this cannot be established, and we need further assumptions.
Consider as in the introduction the following conditions 

\begin{eqnarray} 
&(I)& 
\qu r_{inj} (M,g) = \inf_{x \in M} r_{inj}(x) > 0 , \nonumber \\
&(B_k(M,g))& 
\qu |(\nabla^g)^i R^g| \le C_i, \en 0 \le i \le k ,\nonumber \\
&(B_k(E,\nabla^E))& 
\qu |(\nabla^g)^i R^E| \le C_i, \en 0 \le i \le k .\nonumber
\end{eqnarray}

\no 
It is a well known fact that for any open manifold and given $k, 
0 \le k \le \infty$, there exists a metric $g$ satisfying $(I)$ and 
$(B_k(M,g))$. Moreover, $(I)$ implies completeness of $g$.

\me
\no
{\bf Lemma 2.1.} {\it 
Assume $(M^n,g)$ with (I) and $(B_k)$. Then $C^\infty_c(E)$ is a dense
subset of $W^{p,r}(E)$ and $H^{p,r}(E)$ for $0 \le r \le k+2$.
}

\me
\no
See [5], proposition 1.6 for a proof. \hfill $\Box$

\me
\no
{\bf Lemma 2.2.} {\it
Assume $(M^n,g)$ with (I) and $(B_k)$. Then there exists a
continuous embedding
\[ W^{p,r}(E) \hookrightarrow H^{p,r}(E), \quad 0 \le r \le k+1 . \]
}

\me
\no
{\bf Proof.} According to 2.1, we are done if we could prove
\[ |\Phi|_{H^{p,r}} \le C \cdot |\Phi|_{W^{p,r}} \]
for $0 \le r \le k+1$ and $\Phi \in C^\infty_C(E)$. Perform induction.
For $r=0$, $|\Phi|_{H^{p,0}} = |\Phi|_{W^{p,0}}$. Assume 
$|\Phi|_{H^{p,r}} \le C \cdot |\Phi|_{W^{p,r}}$. Then
\begin{eqnarray} 
  |\Phi|_{H^{p,r+1}} & \le & C \cdot ( |\Phi|_{H^{p,r}} +
  |D \Phi|_{H^{p,r}})  \nonumber \\  
  & \le & C \cdot ( |\Phi|_{W^{p,r}} + |D \Phi|_{W^{p,r}}) . \nonumber
\end{eqnarray}
Let 
$\frac{\partial}{\partial x^i}, i=1, \dots, n$ 
be coordinate vectors fields which are orthonormal in $m \in M$. 
Then with 
$\nabla_i = \nabla_{\frac{\partial}{\partial x^i}}$
\[ |\nabla^s D \Phi|^p_m \le C \cdot \sum_{i_1, \dots ,i_s, j} 
   |\nabla_{i_1} \dots \nabla_{i_s} \frac{\partial}{\partial x^j} \cdot
   \nabla_j \Phi|^p . \]   
Now we apply the Leibniz rule and use the fact that in an atlas of 
normal charts the Christoffel symbols have bounded euclidean 
derivatives up to order $k-1$. This yields
\[ |\nabla^r D \Phi|^p_m \le C \cdot \sum_{i_1, \dots ,i_{r+1}} 
   |\nabla_{i_1} \dots \nabla_{i_{r+1}} \Phi|^p_m \en \mbox{for} \en
   r \le k , \]   
i. e.
\[ |D\Phi|_{W^{p,r}} \le C \cdot |\Phi|_{W^{p,r+1}} \]    
altogether
\[ |\Phi|_{H^{p,r+1}} \le C \cdot |\Phi|_{W^{p,r+1}}. \]
$\hfill \Box $

\me
\no
{\bf Remark.} For $p=2$ this proof is contained in [2]. \hfill $\Box$

\me
\no
{\bf Theorem 2.3.} {\it Assume $(M^n,g)$ with (I) and $(B_k)$ and
$(E,\nabla)$ with $(B_k)$ and $p=2$. Then for $r \le k$
\[ H^{2,r}(E) \equiv H^r(E) \cong W^r(E) \equiv W^{2,r}(E) \]
as equivalent Hilbert spaces.}

\me
\no
{\bf Proof.} According to 2.2., $W^r(E) \subseteq H^r (E)$
continuously. Hence we have to show $H^r(E) \subseteq W^r(E)$
continuously. The latter follows from the local elliptic 
inequality, a uniformly locally finite cover by normal charts of
fixed radius,  trivializations and the existence of 
 elliptic constants. The proof is performed in [2].
$\hfill \Box$

\me
\no
{\bf Remark.} 2.3 holds for $1 < p < \infty$ (cf. [13]). 
$\hfill \Box $

\me
As it is clear from the definition, the spaces $W^{p,k}(E)$ can be
defined for any Riemannian vector bundle $(E, h_E, \nabla^E)$. We
assume this more general case and define additionally
 
\[ {}^{b,s}W(E) := \left\{ \varrho \in C^S(E) \en \Big| \en 
   {}^{b,s}|\varrho| := \sum^s_{i=0} \sup_{x \in M}  
   |\nabla^i \varrho|_x < \infty   \right\} \]          
                                                      
\no 
and in the case of a Clifford bundle
 
\[ {}^{b,s}H(E) := \left\{ \varrho \in C^S(E) \en \Big| \en 
   {}^{b,s,D}|\varrho| := \sum^s_{i=0} \sup_{x \in M}  
   |D^i \varrho|_x < \infty   \right\} . \]          
                                                      
\no
${}^{b,s}W(E)$ is a Banach space and coincides with the completion
of the space of all $\varrho \in \Gamma(E)$ with 
${}^{b,s}|\varrho| < \infty$ with respect to ${}^{b,s}|{\, }|$.

\me 
\no
{\bf Theorem 2.4.} {\it Let $(E,h,\nabla^E)$ be a Riemannian vector
bundle satisfying (I), $(B_k(M^n,g))$, $B_k(E;\nabla))$.

\setcounter{equation}{0}   

\no
{\bf a.} Assume $k \ge r, k \ge 1, r-\frac{n}{p} \ge s-\frac{n}{q},
r \ge s, q \ge p $, then
\be
W^{p,r}(E) \hookrightarrow W^{q,s}(E) 
\ee
continuously.

\no
{\bf b.} If $k \ge 0, r > \frac{n}{p} + s $ then
\be 
W^{p,r}(E) \hookrightarrow {}^{q,s}W(E) 
\ee
continuously.}

\me
\no
We refer to [7] for the proof. $\hfill \Box$

\me
\no
{\bf Corollary 2.5.} {\it Let $E \rightarrow M$ be a Clifford bundle
satisfying (I), $(B_k(M))$, $(B_k(E))$, $k>r>\frac{n}{2}+s$. Then
\be 
H^r(E) \hookrightarrow {}^{b,s}H(E) 
\ee
continuously.}                      

\me
\no
{\bf Proof.} We apply 2.3, (2.2) and obtain
\be
H^r(E) \hookrightarrow {}^{b,s}W(E). 
\ee    
Quite similar as in the proof of 2.2., 
\be
H^r(E) \hookrightarrow {}^{b,s}W(E). 
\ee    
continuously. \hfill $\Box$

\me
A key role for anything in the sequel plays the module structure theorem
for Sobolev spaces. 

\sloppy
\me 
\no
{\bf Theorem 2.6.}
{\it 
Let $(E_i,h_i,D_i) \rightarrow (M^n,g)$
  be vector bundles with $(I)$, $(B_k(M^n,g))$, $(B_k(E_i,\nabla_i))$,
  $i=1,2$. Assume $0\le r\le r_1,r_2\le k$. If $r=0$ assume
  \[
  \left\{ 
    \begin{array}{rcl}
      r-\frac{n}{p} & < & r_1-\frac{n}{p_1} \\
      r-\frac{n}{p} & < & r_2-\frac{n}{p_2} \\
      r-\frac{n}{p} & \le & r_1-\frac{n}{p_1} + r_2-\frac{n}{p_2} \\
      \frac{1}{p} & \le & \frac{1}{p_1} +\frac{1}{p_2}
    \end{array}
  \right\}
  \en \mbox{or} \]
\be
  \left\{ 
    \begin{array}{rcl}
      r-\frac{n}{p} &  \le  & r_1-\frac{n}{p_1} \\
      0  & <  & r_2-\frac{n}{p_2} \\
      \frac{1}{p} & \le  & \frac{1}{p_1}
    \end{array}
  \right\}
  \en \mbox{or} \en
  \left\{ 
    \begin{array}{rcl}
      0  & <  & r_1-\frac{n}{p_1} \\
      r-\frac{n}{p} &  \le  & r_2-\frac{n}{p_2} \\
      \frac{1}{p} & \le  & \frac{1}{p_2}
    \end{array}
  \right\}.
\ee
  If $r>0$ assume $\frac{1}{p}\le \frac{1}{p_1} + \frac{1}{p_2}$ and
\be
  \left\{ 
    \begin{array}{rcl}
      r-\frac{n}{p} & < & r_1-\frac{n}{p_1} \\
      r-\frac{n}{p} & < & r_2-\frac{n}{p_2} \\
      r-\frac{n}{p} & \le & r_1-\frac{n}{p_1} + r_2-\frac{n}{p_2}
    \end{array}
  \right\}
  \en \mbox{or} \en
  \left\{ 
    \begin{array}{rcl}
      r-\frac{n}{p} & \le & r_1-\frac{n}{p_1} \\
      r-\frac{n}{p} & \le & r_2-\frac{n}{p_2} \\
      r-\frac{n}{p} & < & r_1-\frac{n}{p_1} + r_2-\frac{n}{p_2}
    \end{array}
  \right\}.
\ee
Then the tensor product of sections defines a continuous bilinear map
\[ W^{p_1,r_1}(E_1,\nabla_1) \times W^{p_2,r_2}(E_2,\nabla_2) \longrightarrow
   W^{p,r} (E_1 \otimes E_2, \nabla_1 \otimes \nabla_2) . \]
}
   
\me
\no
We refer to [7] for the proof. \hfill $\Box$

\me
Define for $u \in C^0(M), c>0$ 
\[ \overline{u}_c(x) := \frac{1}{vol B_c(x)} 
   \int\limits_{B_c(x)} u(y) d vol_y (g) . \]

\me 
\no
{\bf Lemma 2.7.}
{\it
Let $(M^n,g)$ be complete $Ric(g) \ge k, k \in {\bf R}$.
Then there exists a positive constant $C=C(n,k,R)$, depending only on 
$n,k,R$ such that for any $c \in  ]0,R[$ and any 
$u \in W^{1,1}(M) \cap C^\infty (M)$
\[ \int\limits_M |u-\overline{u}_c| d vol_x (g) \le  
   C \cdot c \cdot \int\limits_M |\nabla u| d vol_x (g) .  \]
}

\me
\no
{\bf Proof.} For $u \in C_c^\infty (M)$ the proof is performed in [10],
p. 31--33. But what is only needed in the proof is 
$\int |u| dx, \int |\nabla u| dx < \infty $
(even only 
$\int |\nabla u| dx < \infty $).
   
\me
The key is the lemma of Buser, 
\[ \int\limits_{B_c(x)} |u-\overline{u}_c| dx \le
   C \cdot c \cdot \int\limits_{B_c(x)} |\nabla u| dx . \]
\hfill $\Box$

\me
\no
{\bf Remark.} Even $u, \nabla u \in C^\infty$ are not necessary, 
$u \in C^1$ is completely sufficient. 

\hfill $\Box$

\me
\no
{\bf Proposition 2.8.}
{\it
Let $(E,h,\nabla) \rightarrow (M^n,g)$ be a Riemannian vector bundle,
$(M^n,g)$ with (I), $(B_0), r>n+1, 0<c<r_{inj}$ and 
$\eta \in W^{1,r}(E)$. Then 
$\overline{|\eta|}_c \in W^{1,0}(M) \equiv L_1(M)$, where
\[ \overline{|\eta|}_c(x) := \frac{1}{vol B_c(x)}
   \int\limits_{B_c(x)} |\eta(y)| dy . \]
}

\no
{\bf Proof.} Set $u(x)=|\eta(x)|$. Then   
$\overline{u_c(x)}=\overline{|\eta|}_c(x)$ and, according to Kato's
inequality, 
\[ \int |\nabla u| dx = \int |\nabla |\eta|| dx \le 
   \int |\nabla \eta| dx < \infty . \]
Hence we obtain from 2.7, $|u| = |\eta| \in L_1, 
|\eta| - \overline{|\eta|_c} \in L_1$,
\be
\overline{|\eta|}_c \in L_1 .
\ee
\hfill $\Box$

\me
\no
{\bf Remark.} For (2.2) is the assumption $(B_0(E))$ superfluous.
Nevertheless, we have in our applications even $(B_k(E))$.
\hfill $\Box$

\me
Finally we recall for clarity and distinctness a fact which will be
very important later. Let $(E,h,\nabla) \rightarrow (M^n,g)$ be a
Riemannian vector bundle with $(I), (B_k(M)), (B_k(E)), k \ge r+1,
r > \frac{n}{p}+1, 0 < c < r_{inj}$. Then the spaces 
$W^{p,r}(E |_{B_c(x)}) = \{ \varrho \Big| \varrho$ distributional section
of $E |_{B_c(x)}$ s. t. $|\varrho|_{p,r} < \infty \}$ are well defined,
$x \in M$ arbitrary. Radial parallel translation of an orthonormal basis
defines an isomorphism

\be
A_x : W^{p,r} (E |_{B_c(x)}) \stackrel{\cong}{\longrightarrow}
W^{p,r}(B_c(0),V^N),
\ee                 

\no
$B_c(0) \subset {\bf R}^n, V^N = {\bf R}^N$ or ${\bf C}^N, N = rk E$.
We conclude from $(B_k(M)), (B_k(E))$, $k \ge r+1$ and [7] that there
exists constants $c_1, C_1$ s. t.

\be
c_1 \cdot |\varrho|_{p,r,B_c(x)} \le  |A_x \varrho|_{p,r,B_c(0)}  
\le C_1 \cdot |\varrho|_{p,r,B_c(x)} ,  
\ee                                   

\no
$c_1, C_1$ independent of $x$. Moreover, if $\varrho \in W^{p,r}(E)$
then $\varrho |_{B_c(x)} \in W^{p,r} (E |_{B_c(x)})$. Similarly, 

\be
c_2 {}^{b,s} |\varrho|_{B_c(x)} \le {}^{b,s} |A_x \varrho|_{B_c(0)}  
\le C_2 \cdot {}^{b,s}|\varrho|_{B_c(x)} ,  
\ee

\no
$c_2, C_2$ independent of $x$. $(B_k(M)), (B_k(E)),0 < c < r_{inj}$
imply that $B_c(x)$ satisfies all required smoothness conditions and we
obtain from the Sobolev embedding theorem, (2.10), (2.11)

\be
\hspace*{0.7cm} 
W^{p,r} (E |_{B_c(x)}) \hookrightarrow {}^{b,1} W(E |_{B_c(x)}) ,
\ee                                 

\be
{}^{b,1}|\varrho|_{B_c(x)} \le C \cdot |\varrho|_{p,r,B_c(x)} ,
\ee

\no
$C$ independent of $x$.

\section{Uniform spaces of metrics and connections}

\setcounter{equation}{0}

Denote by ${\cal M}(I,B_k)$ the set of all metrics $g$ satisfying the
conditions $(I)$ and $(B_k)$. 

Let $1 \le p < \infty, k \ge r \ge \frac{n}{p}+2, \delta > 0$ and set
\begin {eqnarray} 
   V_\delta &=& \Bigg\{ (g,g') \in {\cal M}(I,B_k)^2 \Big| \en  g \en \mbox{and} \en 
   g' \en \mbox{are quasi isometric and} \en
   |g-g'|_{g,p,r} \nonumber \\
   &{}& := \en \left( \int \left( |g-g'|^p_{g,x} +
   \sum^{r-1}_{i=0} | (\nabla^g)^i (\nabla^g-\nabla^{g'}) |^p_{g,x}
   \right) dvol_x(g)
   \right)^\frac{1}{p} < \delta
   \Bigg\} . \nonumber
\end{eqnarray}

\no 
Here $g,g'$ quasi isometric means $C_1 \cdot g \le g' \le C_2 \cdot g$
in the sense of quadratic forms. This is equivalent to 
$ {}^b|g-g'|_g < \infty $ and ${}^b|g-g'|_{g'} < \infty $, 
where for a tensor 
$t$ ${}^b|t|_g = \sup_{x \in M} |t|_{g,x}$.

\me
\no
{\bf Proposition 3.1.} {\it
Assume $p,k,r$ as above. Then ${\cal B} = \{ V_\delta \}_{\delta>0}$
is a basis for a metrizable uniform structure 
${\cal U}^{p,r}({\cal M}(I,B_k))$.
}
  
\no                                                                             
We refer to [6] for a proof. The key to the proof is the module
structure theorem. 

\hfill $\Box$

\me
Let ${\cal M}^p_r(I,B_k) = {\cal M}(I,B_k)$ 
endowed with the  topology.
${\cal M}^{p,r} := \overline{{\cal M}^p_r}$ the completion. If 
$k \ge r > \frac{n}{p}+1$ then ${\cal M}^{p,r}$ 
still consists of $C^1$--metrics,
i.e. does not contain semi definite elements. This has been proved by
Salomonsen in [12].

\me
\no
{\bf Theorem 3.2.} {\it
Let $k \ge r > \frac{n}{p}+2, g \in {\cal M}(I,B_k)$, 
$U^{p,r}(g) = \Big\{ g' \in {\cal M}^{p,r}(I,B_k) \big|$ ${}^b|g-g'|_g < 
\infty, \en
{}^b|g-g'|_{g'} < \infty$ and $|g-g'|_{g,p,r} < \infty \Big\}$
and denote by $comp(g) \subset {\cal M}^{p,r}(I,B_k)$ the component of 
$g$ in ${\cal M}^{p,r}(I,B_k)$. Then 
\be
comp(g) = U^{p,r}(g) ,
\ee
$comp(g)$ is a Banach manifold, for $p=2$ a Hilbert manifold and
${\cal M}^{p,r}(I,B_k)$ has a representation as topological sum
\be
{\cal M}^{p,r}(I,B_k) = \sum_{j \in J} comp(g_j)
\ee
J an uncountable set.
}

\me
\no
The proof is performed in [6]. \hfill $\Box$

\me
\no
{\bf Remarks.}

\no
{\bf 1.} If $M^n$ is compact then $J$ consists of one element.

\no
{\bf 2.} All metrics in the completed space are at least of class 
$C^2$. Hence curvature is well defined. 

\hfill $\Box$

\me
Let $(E,h) \rightarrow (M^n,g)$ be a Clifford bundle without a fixed
connection, $(M^n,g)$ with $(I)$ and $(B_k)$. 

Set

$C_E(B_k) = \Big\{ \nabla \Big|$ is Clifford connection, metric with
respect to $h$ and satisfies 

\hspace{2cm}$(B_k(E,\nabla)) \Big\}$

Assume $(E,h) \rightarrow (M^n,g)$ as above, $k \ge r > \frac{n}{p}+2,
\delta >0$ and set 

\begin{eqnarray} 
  V_\delta &=& \Big\{ (\nabla, \nabla') \in C_E(B_k)^2 \Big| 
  |\nabla-\nabla'|_{\nabla,p,r} \nonumber \\
  &{}& := \en \Big( \int \sum^r_{i=0} |\nabla (\nabla - \nabla')|^p_x
  dvol_x(g) \Big)^{\frac{1}{p}} < \delta \Big\} . \nonumber
\end{eqnarray}

\me
\no
{\bf Proposition 3.3.} {\it
Assume $p,k,r$ as above. Then ${\cal B} = \{V_\delta \}_{\delta>0}$
is a basis for a metrizable  structure ${\cal U}^{p,r}(C_E(B_k))$.
}

\me 
\no
We refer to [4] for a proof. \hfill $\Box$

\me  
Let $C^{p,r}_E(B_k)$ be the completion of $(C_E(B_k), 
{\cal U}^{p,r}(C_E(B_k)))$. If $\nabla, \nabla' \in C^{p,r}_E(B_k)$ then
$\nabla-\nabla'$ is a 1--form $\eta$ with values in ${\cal G}_E$ =
skew endomorphisms satisfying
\be
\eta_x(Y \cdot \Phi) = Y \cdot \eta_x(\Phi) .
\ee
   
As well known, a metric connection $\nabla$ in $E$ induces a connection
$\nabla$ in ${\cal G}_E$. Denote
\begin{eqnarray}
   \Omega^1({\cal G}^{Cl}_E) &:=& \big\{ \eta \in \Omega^1 ({\cal G}_E)
   \Big| \eta \en \mbox{satisfies (3.3)} \en \big\} , \nonumber \\
   \Omega^{1,p}_r({\cal G}^{Cl}_E,\nabla) &:=&
   \Big\{ \eta \in \Omega^1 ({\cal G}^{Cl}_E)
   \Big| \nonumber \\
   &{}& |\eta|_{\nabla,p,r} 
   := \Big( \int \sum^r_{i=0} |\nabla^i \eta|^p_x
   dvol_x(g) \Big)^{\frac{1}{p}} < \infty \Big\} , \nonumber \\
   \Omega^{1,p,r}({\cal G}^{Cl}_E,\nabla) &:=&
   {\overline{\Omega^{1,p}({\cal G}^{Cl}_E,\nabla)}}^{| |_{\nabla,p,r}} . \nonumber
\end{eqnarray}
        
\no
If $(M^n,g)$ satisfies $(I), (B_k)$ then
\be     
\Omega^{1,p,r}({\cal G}^{Cl}_E,\nabla) \en := \en      
{\overline{C^\infty_c ({\cal G}_E)}}^{| |_{\nabla,p,r}} =
\Big\{ \eta \en \mbox{distributional} \en  \Big| |\eta|_{\nabla,p,r} 
< \infty \Big\} ,
\ee     

$ r \le k+2 $. 
        
\me
\no     
{\bf Theorem 3.4.} {\it
Assume $(E,h) \rightarrow (M^n,g), p, k, r$ as above. Denote for
$\nabla \in C_E(B_k)$ by $comp(\nabla) \subset C^{p,r}_E(B_k)$ the component
of $\nabla$ in $C^{p,r}_E(B_k)$. Then
\be
  comp(\nabla) = \nabla + \Omega^{1,p,r,} ({\cal G}^{Cl}_E, \nabla)
\ee     
and $C^{p,r}_E(B_k)$ has a representation as topological sum 
\be     
 C^{p,r}_E(B_k) = \sum_{j \in J} comp(\nabla_j) .
\ee
}        
         
\me      
\no      
The proof is performed in [4]. \hfill $\Box$
         
\me      
\no      
{\bf Remarks.}
         
\no      
{\bf 1.} If $M^n$ is compact then $C^{p,r}_E(B_k) = C^{p,r}_E$ consists
of one component
         
\no      
{\bf 2.} If $\nabla$ is not smooth then one 
sets $\nabla^i=(\nabla_0+(\nabla-\nabla_0))^i$,
$\nabla_0 \in comp{\nabla} \cap C_E(B_k)$, 
and the right hand side makes sense.
         
\no      
{\bf 3.} All connections in the complete space are at least of class $C^2$.
Hence curvature is well defined.
         
\me      
For the sequel, we must sharpen our considerations concerning Sobolev 
spaces. Let $(E,h,\nabla) \rightarrow (M^n,g)$ be a Riemannian vector
bundle. The connection $\nabla$ enters into the definition of the Sobolev
spaces $W^{p,r}$. Hence we should write $W^{p,r}(E,\nabla)$. Now there 
arises the natural question, how do the spaces $W^{p,r}(E,\nabla)$ on 
$\nabla$? We present here one answer. Other considerations are performed
in [4], [3].
         
\me                                                    
\no       
{\bf Proposition 3.5.} {\it
Let $(E,h,\nabla) \rightarrow (M^n,g)$ be a Riemannian vector
bundle with (I), $(B_k)$, $(B_k(E,\nabla))$, $k \ge r > \frac{n}{p}+1,
1 \le p < \infty$. Suppose $\nabla' \in comp(\nabla) \subset C^{p,r}_E(B_k)$,
$\nabla'$ smooth, i. e. $\nabla'=\nabla+\eta, \eta \in 
\Omega^{1,p,r}({\cal G}_E,\nabla) \cap C^\infty $.
Then              
\be               
W^{p,i}(E,\nabla) = W^{p,i}(E,\nabla'), \en 0 \le i \le r,
\ee              
as equivalent Banach spaces.
}                 
  
\me
For the proof we refer to [4], [3]. The proof includes some combinatorial
considerations and essentially uses the module structure theorem. This
is the reason why we assumed $k \ge r > \frac{n}{p}+1$. But this
assumption can be weakened. We only need the validity of the module
structure theorem. \hfill $\Box$

\me
\no
{\bf Remark.} The assumption $\eta$ smooth in superfluous. As we 
mentioned already several times, we can define $W^{p,r}(E,\nabla')$
and prove (3.7) for $\nabla' \in comp(\nabla)$ only. \hfill $\Box$

\me
\no
{\bf Corollary 3.6.} {\it
Suppose $(E,h,\nabla) \rightarrow (M^n,g)$ as above, 
$k \ge r >\frac{n}{p}+2$, $\nabla' = \nabla + \eta$, 
$\eta \in \Omega^{1,p,r}({\cal G}_E,\nabla)$. Then
\be                                    
W^{2p,i}(E,\nabla) = W^{2p,i}(E,\nabla'), \en
0 \le i \le \frac{r}{2} .
\ee
}

\me
\no
{\bf Proof.} $r>\frac{n}{p}+2$ implies 
$r-\frac{n}{p} \ge \frac{r}{2} - \frac{n}{2p} , 2p \ge p , 
r \ge \frac{r}{2}$, i. e. 
\[ \Omega^{1,p,r}({\cal G}_E,\nabla) \subseteq 
   \Omega^{1,2p,\frac{r}{2}}({\cal G}_E,\nabla). \]
Now we apply 3.5 replacing $p \rightarrow 2 p, r \rightarrow
\frac{r}{2}$. \hfill $\Box$

\me
\no
{\bf Corollary 3.7.} {\it
Suppose $(E,h,\nabla) \rightarrow (M^n,g)$ a Clifford bundle with the
conditions above for $p=1$, i. e. $k \ge r > n+2$ , 
$\nabla' \in comp(\nabla) \subset C^{1,r}_E(B_k)$, $\nabla'$ smooth. Then
\be
W^{2,i}(E,\nabla) \equiv W^i(E,\nabla) = W^i(E,\nabla') \equiv W^{2,i}(E,\nabla'),
\en 0 \le i \le \frac{r}{2} .
\ee
}

\me
\no
{\bf Corollary 3.8.} {\it
Assume the hypothesises of 3.7. and write $D=D(\nabla,g), D'=D(\nabla',g)$.
Then
\be
H^i(E,D) = H^i(E,D'), \en 0 \le i \le \frac{r}{2} 
\ee
In particular,
\be
{\cal D}_{D^i} = {\cal D}_{D'^i}, \en 0 \le i \le \frac{r}{2}  
\ee           
where ${\cal D}_{D^i}$ denotes the domain of definition of
$\overline{D^i}$ .
}

\me
\no
{\bf Proof.} (3.11) follows from the result of Chernoff that
$D^i$ is essentially self adjoint on $C^\infty_c(E)$,
${\cal D}_{D^i} = H^i(E,D)$ and (3.10). (3.10) follows from (3.9) and 2.3.
\hfill $\Box$

\me
Finally we make some remarks concerning the essential spectrum of
$D$ and $D^2$. More precisely, we prove that it is an invariant of 
$comp(\nabla)$. We have several distinct proofs for this and present
here a particularly simple one.

We consider Weyl sequences and restrict to orthonormal ones. Denote
by $\sigma_e(D)$ the essential spectrum of $D$. 
$\lambda \in \sigma_e(D)$ if and only if there exists a Weyl 
sequence for $\lambda$, i. e. an orthonormal sequence 
$(\Phi_\nu)_\nu, \Phi_\nu \in {\cal D}_D$; s. t. 
\be
\lim_{\nu \rightarrow \infty} (D-\lambda) \Phi_\nu = 0 .
\ee

\me
\no
{\bf Lemma 3.9.} {\it 
Suppose $\lambda \in \sigma_e(D)$. Then there exists a Weyl sequence
$(\Phi_\nu)_\nu$ for $\lambda$ s. t. for any compact subset
$K \subset M$
\be
\lim_{\nu \rightarrow \infty} |\Phi_\nu|_{L_2(K,E)} = 0 .
\ee
}

\me
\no
This is Lemma 4.29 of [2]. One simply chooses an exhaustion 
$K_1 \subset K_2 \subset \dots , \bigcup K_i = M$, starts with an
arbitrary Weyl sequence $(\Psi_\nu)_\nu$, produces by the Rellich
lemma and a diagonal choice a subsequence $\chi_\nu$ such that
$(\chi_\nu)_\nu$ converges on any $K_i$ in the $L_2$-sense and defines
$\Phi_\nu := (\chi_{2\nu+1}-\chi_{2\nu}) / \sqrt{2}$. $(\Phi_\nu)_\nu$    
has the desired properties. \hfill $\Box$

\me
\no
{\bf Proposition 3.10.} {\it
Suppose $(E,h,\nabla) \rightarrow (M^n,g)$ a Clifford bundle with (I),
$(B_k(M))$, $(B_k(E,\nabla))$, $k \ge r > n+2, n \ge 2$,
$\nabla' \in comp(\nabla) \subset C^{1,r}_E(B_k)$,
$D=D(\nabla,g), D'=D(\nabla',g)$. 
Then
\be
\sigma_e(D) = \sigma_e(D') .
\ee
}

\me
\no
{\bf Proof.}
\[ D'= \sum_i e_i \nabla'_{e_i} = \sum_i e_i \cdot (\nabla_{e_i} + 
   \eta_{e_i}(\cdot)) 
   = D + \eta^{op} ,\]
where the operator $\eta^{op}$ acts as
\[ \eta^{op}(\Phi)|_x = \sum_i e_i \cdot \eta_{e_i}(\Phi)|_x . \]
Then, pointwise, $|\eta^{op}|_x \le C \cdot |\eta|_x$,
$C$ independent of $x$. Given $\varepsilon > 0$, there exists a compact set
$K=K(\varepsilon) \subset M$ such that
\be                                                        
\sup_{x \in M \setminus K} |\eta|_x < \frac{\varepsilon}{C}, \en 
i. e. \en \sup_{x \in M \setminus K} |\eta^{op}|_x < \varepsilon .    
\ee
Assume now $\lambda \in \sigma_e(D), (\Phi_\nu)_\nu$ a Weyl sequence
as in (3.13). According to (3.11), $\Phi_\nu \in {\cal D}_{D'}$.
Then 
\[ (D'-\lambda) \Phi_\nu = (D'_D) \Phi_\nu + (D-\lambda) \Phi_\nu . \]
By assumption, $(D-\lambda)\Phi_\nu \rightarrow 0$. Moreover,
\[ |(D'-D)\Phi_\nu|_{L_2(M,E)} =
   |\eta^{op}\Phi_\nu|_{L_2(M,E)} \le
   C \cdot (|\eta\Phi_\nu|_{L_2(K,E)} + 
   |\eta\Phi_\nu|_{L_2(M \setminus K,E)}) . \]
$|\eta \Phi_\nu|_{L_2(K,E)} \rightarrow 0$ and
\[ C \cdot |\eta \Phi_\nu|_{L_2(M \setminus K,E)} \le
   C \cdot \sup_{x \in M \setminus K} |\eta|_x  \cdot
   |\Phi_\nu|_{L_2(M \setminus K,E)} < \varepsilon . \]
Hence
$(D'-\lambda)\Phi_\nu \rightarrow 0$,
$\lambda \in \sigma_e(D')$,
$\sigma_e(D) \subseteq \sigma_e(D')$.
Exchanging the role of $D,D'$, we obtain 
$\sigma(D') \subseteq \sigma(D)$.
\hfill $\Box$

\section{General heat kernel estimates}
\setcounter{equation}{0}

We collect some standard facts concerning the heat kernel of 
$e^{-tD^2}$. The best references for this are [1], [2].

We consider the self-adjoint closure of $D$ in $L_2(E)=H^0(E)$,
$D=\int\limits^{+ \infty}_{- \infty} \lambda E_\lambda$.

\me
\no
{\bf Lemma 4.1.} {\it
$\{e^{itD}\}_{t \in {\bf R}}$ defines a unitary group on the spaces
$H^r(E)$, for $0 \le h \le r$ holds
\be
|D^h e^{itD} \Psi|_{L_2} =  |e^{itD} D^h \Psi|_{L_2} =
|D^h \Psi|_{L_2} .
\ee
}
\hfill $\Box$

\me
We can extend this action to $H^{-r}(E)$ by means of duality.

\me 
\no
{\bf Lemma 4.2.} {\it
$e^{-tD^2}$ maps $L_2(E) \equiv H^0(E) \rightarrow H^r(E)$
for any $r>0$ and
\be
|e^{-tD^2}|_{L_2 \rightarrow H^r} \le C \cdot t^{- \frac{r}{2}}, \en
t \in ]0, \infty[, \en C=C(r) .
\ee
}

\me
\no
{\bf Proof.} Insert into $e^{-tD^2} = \int e^{-t \lambda^2} dE_\lambda $
the equation
\[ e^{-t \lambda^2} = \frac{1}{\sqrt{4 \pi t}} 
   \int\limits^{+ \infty}_{- \infty} e^{i \lambda s} 
   e^{- \frac{s^2}{4t}} ds \]
and use
\[ \sup |\lambda^r e^{-t \lambda^2}| \le C \cdot t^{- \frac{r}{2}} . \]
\hfill $\Box$

\me
\no
{\bf Corollary 4.3.} {\it
Let $r, s \in {\bf Z}$ be arbitrary. Then for $t>0$
$e^{-tD^2}: H^r(E) \rightarrow H^s(E)$ continuously.
}

\sloppy
\me
\no
{\bf Proof.} This follows from 4.2., duality and the semi group property
of $\{e^{-tD^2}\}_{t>0}$. 
\hfill $\Box$

\me
$e^{-tD^2}$ has a Schwartz kernel 
$W \in \Gamma ({\bf R}_+ \times M \times M, E \Box \hspace{-0.42cm} \times E)$, 
\[ W(t,m,p) = \langle \delta(m), e^{-tD^2} \delta(p) \rangle ,\]
where $\delta (m) \in H^{-r}(E) \otimes E_m$ is the map
$\Psi \in H^r(E) \rightarrow \langle \delta(m), \Psi \rangle = \Psi(m)$,
$r> \frac{n}{2}$. The main result of this section is the fact that for
$t>0, W(t,m,p)$ is a smooth integral kernel in $L_2$ with good decay
properties if we assume bounded geometry.

Denote by $C(m)$ the best local Sobolev constant of the map 
$\Psi \rightarrow \Psi(m), r > \frac{n}{2}$, and by $\sigma(D^2)$
the spectrum.

\me
\no
{\bf Lemma 4.4.} {\it

\no
{\bf a.} $W(t,m,p)$ is for $t>0$ smooth in all variables.

\no
{\bf b.} For any $T>0$ and sufficiently small $\varepsilon > 0$ there 
exists $C > 0$ such that
\be
|W(t,m,p)| \le e^{-(t-\varepsilon) \inf \sigma(D^2)} \cdot
C \cdot C(m) \cdot C(p) \en \mbox{for all} \en t \in ]T, \infty[ .
\ee

\no
{\bf c.} Similar estimates hold for 
$(D^i_m D^j_p W)(t,m,p)$. 
}

\me
\no
{\bf Proof.} 

a. First one shows $W$ is continuous, which follows from
$\langle \delta(m), \cdot \rangle$ continuous in $m$ and 
$e^{-tD^2} \delta(p)$ continuous in $t$ and $p$. Then one applies
elliptic regularity.

b. Write 

$|\langle \delta(m), e^{-tD^2} \delta(p) \rangle| =
|\langle (1+D^2)^{-\frac{r}{2}} \delta(m), 
(1+D^2)^r e^{-tD^2} (1+D^2)^{\frac{r}{2}} 
\delta(p) \rangle|$.

c. Follows similar as b. \hfill $\Box$

\me
\no
{\bf Lemma 4.5.} {\it
For any $\varepsilon>0, T>0, \delta>0$ there exists $C>0$ such that for 
$r>0, m \in M, T>t>0$ holds
\be
\int\limits_{M \setminus B_r(m)} |W(t,m,p)|^2 dp \le
C \cdot C(m) \cdot e^{-\frac{(r-\varepsilon)^2}{(4+\delta)t}} .
\ee
A similar estimate holds for $D^i_m D^j_p W(t,m,p)$.  
}                                      
 
\no
We refer to [2] for the proof.
                             
\me               
\no   
{\bf Lemma 4.6.} {\it
For any $\varepsilon>0, T>0, \delta>0$ there exists $C>0$ such that for 
all $m, p \in M$ with $dist(m,p) > 2 \varepsilon, T>t>0$ holds   
\be
|W(t,m,p)|^2  \le
C \cdot C(m) \cdot C(p) \cdot 
e^{-\frac{(dist(m,p)-\varepsilon)^2}{(4+\delta)t}} .
\ee
A similar estimate holds for $D^i_m D^j_p W(t,m,p)$.  
}
 
\no
We refer to [2] for the proof.
                             
\me               
\no   
{\bf Proposition 4.7.} {\it
Assume $(M^n,g)$ with (I) and $(B_k)$, $(E,\nabla)$ with $(B_k)$, 
$k \ge r > \frac{n}{2}+1$. Then all estimates in 4.4, 4.5, 4.6 hold
with  constants.
}

\me
\no
{\bf Proof.} From the assumptions $H^r(E) \cong W^r(E)$ and
$\sup\limits_m C(m) = C =$ global Sobolev constant for $W^r(E)$. 
\hfill $\Box$

\me
Let $U \subset M$ be precompact, open, $(M^+,g^+)$ closed with 
$U \subset M^+$ isometrically and $E^+ \rightarrow M^+$ a
Clifford bundle with $E^+|_U \cong E|_U$ isometrically. Denote by
$W^+(t,m,p)$ the heat kernel of $e^{-t{D^+}^2}$.
 
\me
\no            
{\bf Lemma 4.8.} {\it  
Assume $\varepsilon>0, T>0, \delta>0$. Then there exists $C>0$ such
that for all $T>t>0, m,p \in U$ with $B_{2\varepsilon}(m),
B_{2\varepsilon}(p) \subset U$ holds
\be 
|W(t,m,p)-W^+(t,m,p)| \le 
C \cdot e^{- \frac{\varepsilon^2}{(4+\delta)t}}
\ee
}  
  
\no
We refer to [2] for the simple proof. \hfill $\Box$   
  
\me
\no            
{\bf Corollary 4.9.} {\it   
$tr W(t,m,m)$ has for $t \rightarrow 0^+$ the same asymptotic expansion
as $tr W^+(t,m,m)$.} 

\hfill $\Box$

\section{Trace class property under variation of the Clifford connection}
\setcounter{equation}{0}  

We come now to the first main result of this paper.

\me
\no
{\bf Theorem 5.1.} {\it
Assume $(E,\nabla) \rightarrow (M^n,g)$, $(M^n,g)$ with (I) and $(B_k)$  
$(E,\nabla)$ with $(B_k)$, $k \ge r > n+2, n \ge 2$, 
$\nabla' \in comp(\nabla) \cap C_E(B_k) \subset C^{1,r}_E(B_k)$,
$D=D(g, \nabla), D'=D'(g, \nabla'), $ generalized Dirac operators.
Then
\[ e^{-tD^2} - e^{-tD'^2}\]
is for $t>0$ trace class operator and its trace norm is uniformly
bounded on compact $t$--intervalls $[a_0, a_1], a_0>0$.
}

\me
\no
{\bf Remark.} The condition 
$\nabla' \in comp(\nabla) \cap C_E(B_k) \subset C^{1,r}_E(B_k)$, i. e. 
$\nabla' \in comp(\nabla)$ and additionally $\nabla'$ smooth and satisfying 
$(B_k)$
can be weakened to $\nabla' \in comp(\nabla) \subset C^{1,r}_E(B_k)$.
The main reason for this is that we can write
$\nabla'=\nabla'_0+(\nabla'-\nabla'_0)$, $\nabla'_0 \in  C_E(B_k)$, 
$|\nabla'-\nabla'_0|_{1,r,\nabla} < \varepsilon$. Then one can reestablish the 
whole Sobolev theory etc. extensively using the module structure theorem.
We refer to the forthcoming paper [8]. \hfill $\Box$

The proof of theorem 5.1 will occupy the remaining part of this section.
We always assume the assumptions of 5.1. According to (3.11),
\[ {\cal D}_D = {\cal D}_{D'}, \qu {\cal D}_{D^2} = {\cal D}_{D'^2} . \]

\me
\no
{\bf Lemma 5.2.} {\it
Assume $t>0$. Then 
\be
e^{-tD^2} -  e^{-tD'^2} = \int\limits^t_0 e^{-sD^2} 
(D'^2-D^2) e^{-(t-s)D'^2} ds .         
\ee                                              
}

\me
\no
{\bf Proof.} (5.1) means at heat kernel level
\be
W(t,m,p) - W'(t,m,p) = - \int\limits^t_0 \int\limits_M (W(s,m,q),
(D^2-D'^2) W'(t-s,q,p))_q dq ds ,
\ee
where $( , )_q$ means the fibrewise scalar product at $q$ and 
$dq=dvol_q(g)$. Hence for (5.1) we have to prove (5.2). (5.2) is an
immediate consequence of Duhamel's principle. Only for completeness,
we present the proof of (5.2), which is the last of the following
7 facts and implications.

\me
{\bf 1.} For $t>0$ is $W(t,m,p) \in L_2(M,E,dp) \cap {\cal D}_D^2$.

\me
{\bf 2.} If $\Phi, \Psi \in {\cal D}_D^2$ then
$\int (D^2 \Phi, \Psi) - (\Phi, D^2 \Psi) dvol = 0$ 
(Greens formula). 
         
\me
{\bf 3.} $
   ((D^2+ \frac{\partial}{\partial \tau}) \Phi(\tau, g) 
   \Psi(t-\tau, q))_q - (\Phi(\tau,g), (D^2+ \frac{\partial}{\partial t})
   \Psi(t-\tau, q))_q = $

  $ = (D^2 \Phi(\tau,q), \Psi(t-\tau,q))_q - (\Phi(\tau,q),
   D^2 \Psi(t-\tau,q))_q + \frac{\partial}{\partial \tau}
   (\Phi(\tau,g), \Psi(t-\tau,q))_q $.

\me
{\bf 4.}
$ \int\limits^\beta_\alpha \int\limits_M 
   ((D^2+\frac{\partial}{\partial \tau}) \Phi(\tau,q), 
   \Psi(t-\tau, q))_q - (\Phi(\tau,q),  
    (D^2+\frac{\partial}{\partial t}) \Psi(t-\tau,q))_q dq d\tau = 
   $

$  = \int\limits_M [(\Phi(\beta,q), \Psi(t-\beta,q))_q -
   (\Phi(\alpha,q), \Psi(t-\alpha,q))_q ] dq$.

\me 
{\bf 5.}
$ \Phi(t,q) = W(t,m,q), \Psi(t,q) = W'(t,q,p) $ yields 

$  - \int\limits^\beta_\alpha \int\limits_M (W(\tau,m,q), 
   (D^2+\frac{\partial}{\partial t}) W'(t-\tau,q,p))_q dq d\tau = $

$   = \int\limits_M [(W(\beta,m,q), W'(t-\beta,q,p))_q -
   (W(\alpha,m,q), W'(t-\alpha,q,p))_q] dq$  . 

\me
{\bf 6.}
Performing $\alpha \rightarrow 0^+, \beta \rightarrow t^-$ in 5. yields

$ - \int\limits^t_0 \int\limits_M (W(s,m,q),  
   (D^2+\frac{\partial}{\partial t}) W'(t-s,q,p))_q dq ds =
    W(t,m,p) - W'(t,m,p)$ . 

\me
{\bf 7.}
Finally, using $D^2 + \frac{\partial}{\partial t} = 
D^2 - D'^2 + D'^2 + \frac{\partial}{\partial t} $ and
$(D'^2+\frac{\partial}{\partial t}) W' = 0$ we obtain

$ W(t,m,p) - W'(t,m,p) = - \int\limits^t_0 \int\limits_M (W(s,m,q),
(D^2-D'^2) W'(t-s,q,p))_q dq ds $

which is (5.2). \hfill $\Box$

\en
If we write $D'^2-D^2 = D'(D'-D)+(D'-D)D$ then
\begin{eqnarray}
   e^{-tD^2} - e^{-tD'^2} &=& \int\limits^t_0 e^{-sD^2} D'(D'-D)  
   e^{-(t-s)D'^2} ds \en + \nonumber \\                                       
   &+& \int\limits^t_0 e^{-sD^2} (D'-D) D e^{-(t-s)D'^2} ds 
   \en  = \nonumber \\
   &=& \int\limits^t_0 e^{-sD^2} D'\eta e^{-(t-s)D'^2} ds \en + \nonumber \\
   &+& \int\limits^t_0 e^{-sD^2} \eta D e^{-(t-s)D'^2} ds , \nonumber
\end{eqnarray}
                                                        
\no
where $\eta = \eta^{op}$ in the sense of section 3, 
$\eta^{op}(\Psi)|_x = \sum^n_{i=1} e_i \cdot \eta_{e_i}(\Psi)$ and
$|\eta^{op}|_x \le C \cdot |\eta|_x$, $C$ independent of $x$.
We split 
$\int\limits^t_0 = \int\limits^{\frac{t}{2}}_0 + 
\int\limits^t_{\frac{t}{2}}$,

\def\theequation{$I_1$} 
\be   
   e^{-tD^2} - e^{-tD'^2} \en = \en
   \int\limits^{\frac{t}{2}}_0 e^{-sD^2} D'\eta e^{-(t-s)D'^2} ds \en + 
\ee 
\def\theequation{$I_2$}
\be
   \hspace*{2.7cm}
   + \en \int\limits^{\frac{t}{2}}_0 e^{-sD^2} 
   \eta D e^{-(t-s)D'^2} ds \en + 
\ee
\def\theequation{$I_3$}
\be
   \hspace*{2.8cm}
   + \en \int\limits^t_{\frac{t}{2}} e^{-sD^2} 
   D'\eta e^{-(t-s)D'^2} ds \en + 
\ee
\def\theequation{$I_4$} 
\be
   \hspace*{2.2cm}
   + \en \int\limits^t_{\frac{t}{2}} e^{-sD^2} 
   \eta D e^{-(t-s)D'^2} ds . 
\ee

\def\theequation{\thesection.\arabic{equation}}  
\setcounter{equation}{2}

We want to show that each integral $(I_1)-(I_4)$ is a product of 
Hilbert-Schmidt operators and to estimate its Hilbert-Schmidt norm.
Consider the integrands of $(I_3)$ resp. $(I_4)$. Applying a 
Leibniz type rule in 
\[ e^{-sD^2} D'\eta e^{-(t-s)D'^2} , \]
we have to estimate 
\be
\Big( e^{-sD^2} \nabla' \eta \Big) \circ \Big( e^{-(t-s)D'^2}  \Big)
\ee 
and
\be
\Big( e^{-sD^2} \eta \Big) \circ \Big( D' e^{-(t-s)D'^2}  \Big)    
\ee                                                               
Similarly for $(I_4)$ 
\begin{eqnarray}
   e^{-sD^2} \eta D e^{-(t-s)D'^2} &=& 
   e^{-sD^2} \eta ((D-D')+D') e^{-(t-s)D'^2} = \nonumber \\
   &=& \Big( - e^{-sD^2} \eta^2 \Big) \circ \Big( e^{-(t-s)D'^2}  \Big)   
   \en + \\
   &+& \Big( e^{-sD^2} \eta \Big) \circ \Big( D' e^{-(t-s)D'^2}  \Big) . 
\end{eqnarray}

Then according to 4.2,
\be
| e^{-(t-s)D'^2} |_{L_2 \rightarrow H^1} \le C \cdot (t-s)^{-\frac{1}{2}}
\ee
and
\be                                   
| D' e^{-(t-s)D'^2} |_{L_2 \rightarrow L_2} \le 
| D' |_{H^1 \rightarrow L_2} | e^{-(t-s)D'^2} |_{L_2 \rightarrow H^1} \le  
C' \cdot (t-s)^{\frac{1}{2}} .
\ee  
(5.7) and (5.8) estimate the right hand factors in (5.3) -- (5.6). Start
now with the left hand factor in (5.6), $e^{-sD^2} \eta$ and write
\be                                                            
e^{-sD^2} \eta = e^{-\frac{s}{2}D^2} \circ e^{-\frac{s}{2}D^2} \eta =
(e^{-\frac{s}{2}D^2} \circ f) \circ (f^{-1} \circ 
e^{-\frac{s}{2}D^2} \eta) .
\ee
Here $f$ shall be a scalar function which acts by multiplication. The main
point is the right choice of $f$. $e^{-\frac{s}{2}D^2}f$ has the    
integral kernel                                                     
\be                                                               
W(\frac{s}{2},m,p) f(p)                              
\ee                        
and $f^{-1} e^{-\frac{s}{2}D^2}$ has the kernel
\be                                                               
f^{-1}(m) W(\frac{s}{2},m,p) \eta(p)                              
\ee                       
We have to make a choice such that (5.10), (5.11) are square integrable
over $M \times M$ and that their $L_2$--norm is uniformly bounded.

We decompose the $L_2$--norm of (5.10) as
\begin{eqnarray}
  &{}& \int\limits_M \int\limits_M | W(\frac{s}{2},m,p) |^2 |f(m)|^2 dm dp = 
  \nonumber \\
  &{}& \int\limits_M \int\limits_{dist(m,p) \ge c} | W(\frac{s}{2},m,p) |^2 
  |f(m)|^2 dp dm + \\
  &{}& \int\limits_M \int\limits_{dist(m,p) < c} | W(\frac{s}{2},m,p) |^2 
  |f(m)|^2 dp dm 
\end{eqnarray}

We obtain from 4.4 for $s \in ]\frac{t}{2},t[$
\[ (5.13) \le \int\limits_M C_1 |f(m)|^2 vol B_c(m) dm \le
   C_2 \int\limits_M |f(m)|^2 dm \]
and from 4.5
\[ \int\limits_M \int\limits_{dist(m,p) \ge c} | W(\frac{s}{2},m,p) |^2  
   |f(m)|^2 dp dm  \le 
   \int\limits_M C_1 e^{- \frac{(r-\varepsilon)^2}{4+\delta} \frac{2}{s}} 
   |f(m)|^2 dm \le \]
\be   
   \le C_1 \cdot e^{- \frac{(c-\varepsilon)^2}{4+\delta} \frac{2}{s}} 
   \int\limits_M |f(m)|^2 dm, \quad c > \varepsilon .
\ee
                  
Hence the  estimate of $\int\limits_M \int\limits_M
| W(\frac{s}{2},m,p) |^2 |f(m)|^2 dp dm    $   
for $s \in [\frac{t}{2},t]$ is done if 
\[ \int\limits_M |f(m)|^2 dm < \infty . \]

For (5.11) we have to estimate 
\be
\int\limits_M \int\limits_M |f(m)|^{-2} 
|(W(\frac{s}{2},m,p), \eta^{op}(p) \cdot)_p|^2 dp dm
\ee          
We recall a simple fact in Hilbert spaces. Let $X$ be a Hilbert space, 
$x \in X, x \ne 0$. Then $|x|= \sup\limits_{|y|=1} | \langle x,y \rangle |$,
\be
|x|^2 = \Big( \sup_{|y|=1} | \langle x,y \rangle | \Big)^2  .
\ee
This follows from $|\langle x,y \rangle| \le |x| \cdot |y|$
and equality for $y = \frac{x}{|x|}$. We apply this to $E \rightarrow M$,
$X = L_2 (M,E,dp)$, $x=x(m)(p)=(W(t,m,p),  \eta^{op}(p) \cdot)_p
= \eta^{op}(p) \circ W(t,m,p)$ and have to estimate
\be
\sup_{
\begin{array}{c}
\Phi \in C^\infty_c(E) \\
|\Phi|_{L_2} = 1
\end{array}
}      
N(\Phi) = 
\sup_{
\begin{array}{c}
\Phi \in C^\infty_c(E) \\
|\Phi|_{L_2} = 1
\end{array}
} 
|\langle \delta(m), e^{-tD^2} \eta ^{op} \Phi \rangle|
\ee
According to 4.5, for $t>0$
\be
W(t,m,\cdot) \in H^{\frac{r}{2}}(E), \quad 
|W(t,m,\cdot)|_{H^{\frac{r}{2}}} \le C_1(t) .                                     
\ee
Hence we can restrict in (5.17) to
\be
\sup_{
\begin{array}{c}
\Phi \in C^\infty_c(E) \\
|\Phi|_{L_2} = 1 \\
|\Phi|_{H^{\frac{r}{2}}} \le C_1
\end{array}
}      
N(\Phi) 
\ee
In the sequel we estimate (5.19). For doing this, we recall some simple
facts concerning the wave equation
\be
\frac{\partial \Phi_s}{\partial s} = i D \Phi_s,
\quad \Phi_0 = \Phi,
\quad \Phi \in C^1 \en 
\mbox{with compact support.}
\ee
It is well known that (5.20) has a unique solution $\Phi_s$ which ist given
by
\be
\Phi_s = e^{i s D} \Phi
\ee
and
\be
\mbox{supp} \en \Phi_s \subset U_{|s|} \en( \mbox{supp} \en \Phi)
\ee                           
$U_{|s|} = |s|$ -- neighborhood. Moreover,
\be
|\Phi_s|_{L_2} =  |\Phi|_{L_2}, \quad
|\Phi_s|_{H^{\frac{r}{2}}} = |\Phi|_{H^{\frac{r}{2}}}       
\ee
We fix a uniformly locally finite cover
${\cal U} = \{ U_\nu \}_\nu = \{B_d(x_\nu)\}_\nu$
by normal charts of radius $d < r_{inj}(M,g)$ and associated 
decomposition of unity $\{ \varphi_\nu \}_\nu$ satisfying
\be
|\nabla^i \varphi_\nu| \le C \en \mbox{for all} \en \nu, \en
0 \le i \le k+2
\ee 

\setcounter{equation}{25} 
Write
\bea
   N(\Phi) & = & | \langle \delta(m) , e^{r-tD^2} \eta^{op} \Phi \rangle | 
   \nonumber \\
   & = & \frac{1}{\sqrt{4 \pi t}} \en \Big| \langle \delta(m),
   \int\limits^{+ \infty}_{- \infty} e^{\frac{-s^2}{4t}} e^{isD}
   (\eta^{op} \Phi) ds \rangle \Big| \en = \nonumber \\
   & = & \frac{1}{\sqrt{4 \pi t}} \en \Big| 
   \int\limits^{+ \infty}_{- \infty} e^{\frac{-s^2}{4t}} (e^{isD}
   \eta^{op} \Phi) (m) ds \Big| .
\eea
             
We decompose
\be
 \eta^{op} = \sum\limits_\nu \varphi_\nu \eta^{op} \Phi .
\ee  
(5.27) is a locally finite sum, (5.20) linear. Hence
\be
 (\eta^{op})_s = \sum\limits_\nu (\varphi_\nu \eta^{op} \Phi)_s .
\ee  
Denote as above
\[ | \en |_{p,i} \equiv | \en |_{W^{p,i}} , \]  
in particular
\be 
  | \en |_{2,i} \equiv | \en |_{W^{2,i}} \sim | \en |_{H^i} , \qu i \le k .  
\ee
Then we obtain from (5.23), (5.28), (2.1)

\[ |(\varphi_\nu \eta^{op} \Phi)_s |_{H^{\frac{r}{2}}}    
   = |\varphi_\nu \eta^{op} \Phi |_{H^{\frac{r}{2}}}     
   \le C_2 |\varphi_\nu \eta^{op} \Phi|_{2, \frac{r}{2}} \le \]

\be                                               
  \le C_3 | \eta^{op} \Phi|_{2, \frac{r}{2}, U_\nu}     
  \le C_4 | \eta |_{2, \frac{r}{2}, U_\nu}                                                           
  \le C_5 | \eta |_{1, r-1, U_\nu}                
\ee

\me
\no
since $r-1-\frac{n}{1} \ge \frac{r}{2} - \frac{n}{2}, 
r-1 \ge \frac{r}{2} , 2 \ge 1$ for $r > n+2$ and 
$|\Phi|_{H^{\frac{r}{2}}} \le C_1$. This yields together with (2.3), (2.13)
the estimate
\bea   
  |(\eta^{op} \Phi)_s (m)| & \le & C_6 \cdot
   \sum\limits_
   {\begin{array}{c}
      \nu \\
      m \in U_{|s|}(U_\nu)
    \end{array}}      
   |(\varphi_\nu \eta^{op} \Phi)_s|_{2,\frac{r}{2}} \en \le  \nonumber \\
   & \le & C_7 \cdot
   \sum\limits_
   {\begin{array}{c}
      \nu \\          
      m \in U_{|s|}(U_\nu) 
    \end{array}}        
   | \eta |_{1, r-1, U_\nu}  
   \le C_8 \cdot |\eta|_{1, r-1, B_{2d+|s|}(m)} \en = \nonumber \\
  & = & C_8 \cdot vol( B_{2d+|s|}(m)) \cdot
  \left( \frac{1}{vol B_{2d+|s|}(m)} \cdot
  |\eta|_{1, r-1, B_{2d+|s|}(m)} \right) . \nonumber \\
       { }   
\eea 
There exist constants $A$ and $B$, independent of $m$ s. t.
\[ vol( B_{2d+|s|}(m)) \le A \cdot e^{B|s|} . \]
Write 
\be
  e^{- \frac{s^2}{4t}} \cdot  vol( B_{2d+|s|}(m)) \le     
  C_9 \cdot e^{-\frac{9}{10}\frac{s^2}{4t}} ,      
\ee                                                  
thus obtaining
\[  N(\Phi) \le C_{10} \int\limits^\infty_0 
   e^{-\frac{9}{10}\frac{s^2}{4t}}    
   \left( \frac{1}{vol B_{2d+|s|}(m)} \cdot   
   |\eta|_{1, r-1, B_{2d+|s|}(m)} \right) ds . \] 
Now we apply (2.7) with $R=3d+s$ and infer
\bea
   & {} & \int\limits_M  \frac{1}{vol B_{2d+|s|}(m)} \cdot
   |\eta|_{1, r-1, B_{2d+|s|}(m)} dm \en \le \nonumber \\
   & \le & |\eta|_{1,r-1} + C (3d+s) \cdot (2d+s) 
   |\nabla \eta|_{1,r-1} \en \le \nonumber \\
   & \le & |\eta|_{1,r-1} + C (3d+s) \cdot (2d+s) |\eta|_{1,r} .
\eea
$C (3d+s)$ depends on $3d+s$ at most linearly exponentially, i. e.
\[ C (3d+s) \le A_1 e^{B_1 (3d+s)} . \]
This implies
\bea
   \int\limits^\infty_0 
   e^{-\frac{9}{10}\frac{s^2}{4t}}    
   \int\limits_M        
   \frac{1}{vol B_{2d+|s|}(m)} \cdot   
   |\eta|_{1, r-1, B_{2d+|s|}(m)} dm ds \en \le \nonumber \\
   \le \int\limits^\infty_0 
   e^{-\frac{9}{10}\frac{s^2}{4t}} 
   (|\eta|_{1,r-1} +C (3d+s) \cdot (2d+s)
   |\eta|_{1,r}) ds < \infty .
\eea
The function
\bea
   &{}& {\bf R}_+ \times M \rightarrow {\bf R},  \nonumber  \\
   &{}& (s,m) \rightarrow e^{-\frac{9}{10}\frac{s^2}{4t}}
   \left( \frac{1}{vol B_{2d+|s|}(m)} \cdot   
   |\eta|_{1, r-1, B_{2d+|s|}(m)} \right) \nonumber
\eea
is measurable, nonnegative, the integrals (5.33), (5.34) exist, hence
according to the principle of Tonelli, this function is $1$--summable,
the Fubini theorem is applicable and 
\[ \tilde \eta (m) := C_{10} \cdot  \int\limits^\infty_0 
   e^{-\frac{9}{10}\frac{s^2}{4t}}    
   \left( \frac{1}{vol B_{2d+|s|}(m)} \cdot   
   |\eta|_{1, r-1, B_{2d+|s|}(m)} \right) ds \]
is (for $\eta \not\equiv 0)$ everywhere $\not= 0$ and $1$--summable.
We proved
\be
  \int |(W(t,m,p), \eta^{op} \cdot )_p|^2 dp \le \tilde \eta (m)^2 .
\ee
Now we set
\be
  f(m) = (\tilde \eta (m))^{\frac{1}{2}}
\ee
and infer $f(m) \not= 0$ everywhere, $f \in L_2$ and 
\bea 
   &{}& \int\limits_M \int\limits_M f(m)^{-2} 
   |(W(\frac{t}{2},m,p), \eta^{op} \cdot)_p|^2 dp dm \en \le \nonumber \\
   & \le & \int\limits_M \frac{1}{\tilde \eta (m)}  \tilde \eta (m)^2 dm  
  = \int\limits_M  \tilde \eta (m) dm  \le  C_{11} (t) \cdot \sqrt{t} , 
\eea 
since       
$\int\limits^{+\infty}_{-\infty} e^{-ax^2} dx 
= \frac{\sqrt{\pi}}{\sqrt{a}} $,
$\int\limits^{\infty}_0 e^{-ax^2} dx 
= \frac{\sqrt{\pi}}{2 \sqrt{a}} $ ,                                  
in particular for $t \in [a_0,a_1], a_0 >0,$
\be
   \sup_{t \in [a_0,a_1]} C_{11}(t) \cdot \sqrt{t}
   = C_{12} (a_0,a_1) < \infty .
\ee
As shown above, the integrals (5.12), (5.13) can be estimated by constants
$C_{13}(a_0,a_1)$, $C_{14}(a_0,a_1)$. Finally we use for $A$ of trace class,
$B$ bounded. 
\be
   |A \cdot B|_1 \le |A|_1 \cdot |B|_{op}
\ee
and obtain 
\bea
   & {} & \Big| \int\limits^t_{\frac{t}{2}}  (e^{-sD^2} \eta 
   D'e^{-(t-s)D'^2}) ds  \Big|_1  \en \le  \nonumber\\ 
   & \le & \sup_{s \in [\frac{t}{2},t]} |e^{-sD^2} \eta|_1 
   \int\limits^t_{\frac{t}{2}} 
   |D'e^{-(t-s)D'^2}|_{L_2 \rightarrow L_2} ds \en \le \nonumber \\
   & \le &          
   C_{12} (\frac{t}{2},t) (C_{13} (\frac{t}{2},t) + C_{14} (\frac{t}{2},t) )
   \cdot C \cdot \sqrt{t} = C_{15} (\frac{t}{2},t), 
\eea       
i. e. the operator (5.6) is of trace class and for
$t \in [a_0, a_1], a_0 > 0$   , its trace norm is uniformly bounded.
                
We proceed with the expression (5.5). The only distinction here is the 
appearence of $(\eta^{op})^2$ instead of $\eta^{op}$. Then we estimate
as in (5.6), replacing $\eta^{op}$ by $(\eta^{op})^2$. The estimates
become even better. Or we write
\be
   e^{-sD^2} (\eta^{op})^2 e^{-(t-s)D'^2} = 
   \Big( e^{-sD^2} \eta^{op} \Big) \circ \Big( \eta^{op} e^{-(t-s)D'^2} \Big)
\ee
Here $\eta^{op}$ acts in $\eta^{op} e^{-(t-s)D'^2}$ as a bounded operator
according to the module structure theorem. Hence $(I_4)$ is done.
$(I_3)$ can be settled in exactly the same manner. In $(I_2)$ and 
$(I_1)$ we change the role of the factors. Decompose the integrand
of $(I_2)$ as 
\bea
   &{}& e^{-sD^2} \eta D e^{-(t-s)D'^2} \en = \nonumber \\
   &=& - e^{-sD^2} \circ \eta \circ 
   (\eta e^{-(t-s) \frac{D'^2}{2}} \cdot f^{-1}) \circ
   (f \cdot e^{-(t-s) \frac{D'^2}{2}}) \en +  \\
   &+& e^{-sD^2} \circ (\eta \cdot 
   e^{-(t-s) \frac{D'^2}{2}} \cdot f^{-1}) \circ
   (f \cdot D' e^{-(t-s) \frac{D'^2}{2}}) .
\eea  
According to (5.7), (5.8) and the module structure or embedding theorem,
\be
  e^{-sD^2}, \eta^{op}, e^{-sD^2} \circ \eta^{op} \en 
  \mbox{are bounded for} \en s \le \frac{t}{2} .
\ee
The terms
\bea
  &{}& \eta^{op} \circ e^{-(t-s) \frac{D'^2}{2}} \cdot f^{-1} , \\
  &{}& f \cdot e^{-(t-s) \frac{D'^2}{2}}
\eea
and
\be
  f \cdot D' e^{-(t-s) \frac{D'^2}{2}}
\ee
in $s \in [0, \frac{t}{2}]$ can be estimated as 
\[ f^{-1} \cdot e^{-\frac{s}{2}D^2} \eta^{op} \]
and
\[ e^{-\frac{s}{2} D^2} f \]
for $s \in [\frac{t}{2},t]$ .

Finally the integrand of $(I_1)$ can be written
\bea
   &{}& e^{-sD^2} D' \eta e^{-(t-s)D'^2} \en = \nonumber \\
   &=& e^{-sD^2} \circ ((\nabla'^{op} 
   \eta^{op}) e^{-(t-s) \frac{D'^2}{2}} \cdot f^{-1}) \circ
   (f \cdot e^{-(t-s) \frac{D'^2}{2}}) \en + \\
   &+& e^{-sD^2} \circ (\eta^{op} \cdot 
   e^{-(t-s) \frac{D'^2}{2}} \cdot f^{-1}) \circ
   (f \cdot D' e^{-(t-s) \frac{D'^2}{2}}) .
\eea  
The two right terms of (5.48), (5.49) can be estimated as (5.45) -- (5.47).

This finishes the proof of theorem 5.1. \hfill $\Box$
  
\me
For our applications in section 7 we need still the trace class property
of 
\[ D e^{-tD^2} - D' e^{-tD'^2} = e^{-tD^2} D - e^{-tD'^2} D' . \]
Consider 
\bea
   &{}& e^{-tD^2} D - e^{-tD'^2} D' \en = \en e^{-tD^2} (D-D') + 
   (e^{-tD^2}-e^{-tD'^2}) D'\en = \nonumber \\
   &=& e^{-tD^2} (D-D') +
   \int\limits^t_0 e^{-sD^2} D' (D'-D) D' e^{-(t-s)D'^2} ds + 
   \int\limits^t_0 e^{-sD^2} (D'-D) D D' e^{-(t-s)D'^2} ds. \nonumber
\eea
Now
\[ e^{-tD^2}(D-D') = - e^{-tD^2} \eta^{op} = 
   - \Big( e^{-t \frac{D^2}{2}} f \Big) \circ 
   \Big( f^{-1} e^{-t \frac{D^2}{2}} \eta^{op} \Big) , \]                                
$f$ as in (5.36), and we are done,
\be
  |e^{-tD^2} (D-D')|_1  \le C \cdot \sqrt{t} .
\ee
Decompose
\[ \int\limits^t_0 e^{-sD^2} D'(D'-D) D' e^{-(t-s)D'^2} ds =
   \int\limits^{\frac{t}{2}}_0 + \int\limits^t_{\frac{t}{2}} . \]
The estimate of $\int\limits^t_{\frac{t}{2}}$ amounts to that of
\bea
   &{}& \Big( e^{-\frac{s}{2} \frac{D^2}{2}} \cdot f \Big) \circ 
   \Big( f^{-1} e^{-\frac{s}{2} D^2} \circ \nabla' \eta \Big) , \\    
   &{}& D' e^{-(t-s)D'^2} , \\    
   &{}& \Big( e^{-\frac{s}{2} \frac{D^2}{2}} f \Big) \circ 
   \Big( f^{-1} e^{-s \frac{D^2}{2}} \nabla' \eta \Big) , \\   
   &{}& D'^2 e^{-(t-s)D'^2}
\eea
(5.51) can be estimated as (5.9), assuming 
$|\eta|_{1,r+1} < \infty$ , $k \ge r > n+2$, 
(5.52) as (5.7), (5.53) as (5.9). A small difficulty arises with
(5.54) since 
\be
   \big| D'^2 e^{-(t-s)D'^2} \big|_{L_2 \rightarrow H_1} \le
    C  \big( (t-s)^{-\frac{1}{2}} \big)^2 .
\ee
But, considering (5.37), 
we see that
$|$(5.53)$|_1$  generates a factor
$(t-s)^{\frac{1}{2}}$ and we obtain $(t-s)^{-\frac{1}{2}}$ for
integration which doesn't cause any trouble. Quite similar we handle
\bea
  &{}& \int\limits^t_0 e^{-sD^2} (D-D')D D' e^{-(t-s)D'^2} ds 
  \en = \nonumber \\
  &=& - \int\limits^t_0 e^{-sD^2} (D-D')^2 D' e^{-(t-s)D'^2} ds
  +  e^{-sD^2} (D-D') D'^2 e^{-(t-s)D'^2} ds \nonumber \\
  &=& - \int\limits^t_0 (e^{-sD^2} \eta^2) D' e^{-(t-s)D'D^2} ds 
  + \int\limits^t_0 e^{-sD^2} \eta D' e^{-(t-s)D'D^2} ds . \nonumber
\eea
We decompose $\int\limits^t_0 = \int\limits^{\frac{t}{2}}_0 +
\int\limits^t_{\frac{t}{2}}$ and proceed as in (5.51) -- (5.55).
Hence we proved

\en
\no
{\bf Theorem 5.3.} {\it
Assume $(E,\nabla) \rightarrow (M^n,g)$ with $(I), (B_k)$, 
$(E,\nabla)$ with $(B_k)$ $k \ge r > n+3$, $n \ge 2$,
$\nabla' \in comp(\nabla) \cap {\cal C}_E(B_k) \subset 
{\cal C}^{1,r}_E (B_k)$, $D=D(g,\nabla)$ , $D'=D(g,\nabla')$
generalized Dirac operators. Then
\[ e^{-tD^2} - e^{-tD'^2} \]
and
\[ D e^{-tD^2} - D' e^{-tD'^2} \]
are for $t>0$ trace class operators and their trace norm is
uniformly bounded on compact t--intervalls $[a_0,a_1], a_0>0$ .}
\hfill $\Box$

\section{Trace class property for additional variation of the metric}
\setcounter{equation}{0}    

As we know from the definition, $D=D(g,\nabla)=D(g,E,\nabla)$. In section 5 we
considered $D'=D(g,E,\nabla')$. More general, we should consider
$D'=D(g',\nabla')$. But at the first glance, this does not make sense.
Change of $g$ changes the Clifford algebra $Cl_m$, we have now 
$Cl(T_m, g'_m)$ and hence have to consider modules of $Cl(T_m,g'_m)$.
A Clifford bundle associated to $g'$ must consist fibrewise of such 
modules, we arrive at a new bundle $E'$. $E'$ can have a new fibre
metric $h'$. Ne\-ver\-the\-less, locally $E$ and $E'$ are isomorphic. 
Motivated by the consideration that the metric parameters 
$g \rightarrow g'$, $h \rightarrow h'$ move smoothly, we assume that
$E \rightarrow E'$ moves smoothly, $E \cong E'$ as smooth vector bundle.
Hence we indentify $E$ and $E'$, keeping in mind that the fibres 
$E_m$ have different module structures over different algebras. Such
a module structure is given by a section $\cdot$ of 
$\Gamma (Hom(TM \times E, E))$. Endowing $TM$ with $g, \nabla^g $, $E$
with $h, \nabla, Hom(TM \times E, E)$ becomes a Riemannian vector bundle.
Hence
\[ W^{p,r} ( Hom (TM \times E, E), g, h, \nabla) \]  
is well defined. Assuming $g,g',h,h',\nabla, \nabla'$ such that
\[ W^{1,r} ( Hom (TM \times E, E), g, h, \nabla) \cong
   W^{1,r} ( Hom (TM \times E, E), g', h', \nabla') \]                                  
then the condition    

\def\theequation{Clm}  
\be     
   \cdot - \cdot' \in               
   W^{1,r} (Hom(TM \otimes E, E), g',h', \nabla') 
\ee                                       
\def\theequation{\thesection.\arabic{equation}}      
\setcounter{equation}{0}

\no
makes sense.

We make in this section the following

\en
\no
{\bf Assumptions.} $(E,h,\nabla) \rightarrow (M^n,g), 
(E,h',\nabla') \rightarrow (M^n,g')$ Clifford bundles with $(I)$, $(B_k(M))$,
$(B_k(E))$, $k \ge r > n+2$          
\bea
  &{}& g' \in comp(g) \cap {\cal M}(I,B_k) \subset {\cal M}^{1,r}(I,B_k). \\
  &{}& h \en \mbox{and} \en h' \en \mbox{quasi isometric and} \en 
  |h-h'|_{g,h,\nabla,1,r} < \infty, 
  |h-h'|_{g',h',\nabla',1,r} < \infty. \\
  &{}& |\nabla - \nabla'|_{g,h,\nabla,1,r} < \infty . \\  
  &{}& |\nabla - \nabla'|_{g',h',\nabla',1,r} < \infty .                                       
\eea 
and 

\def\theequation{Clm}  
\be     
   \cdot - \cdot' \in               
   W^{1,r} (Hom(TM \otimes E, E), g',h', \nabla') 
\ee                                       
\def\theequation{\thesection.\arabic{equation}}   
\setcounter{equation}{4}   
  
\no      
Here we understand $\nabla - \nabla'$ as a 1--form with values in End $E$.
(6.3) means
\be
  \nabla - \nabla' \in \Omega^{1,1,r} (\mbox{End} \, E, g, h, \nabla^g, \nabla)
\ee
and
\be
  |\nabla - \nabla'|_{g,h,\nabla',1,r} = \int\limits_M
  \sum^r_{i=0} |\nabla^i(\nabla-\nabla')|_{g,x} dvol_x(g).
\ee 
The main result of this section shall be formulated as follows.

\en
\no
{\bf Theorem 6.1.} {\it Let $(E,h,\nabla) \rightarrow (M^n,g)$, 
$(E,h',\nabla') \rightarrow (M^n,g')$ be Clifford bundles with $(I), (B_k(M)),
(B_k(E))$, $k\ge r>n+2$, and (6.1)--(6.4), (Clm). Let $D=D(g,h,\nabla)$, 
$D'=D(g',h',\nabla')$. Then 
\[ e^{-tD^2} - e^{-tD'^2} \]
is of trace class and the trace norm is uniformly bounded on compact
$t$--inter\-valls $[a_0,a_1], a_0>0 $.
}

\en
\no
{\bf Remarks.}

\no 
{\bf 1.} We shall see below that $e^{-tD^2}$ and $e^{-tD'^2}$ act between
the same spaces.

\no
{\bf 2.} For $g=g', h=h'$ we obtain back theorem 5.1. \hfill $\Box$

\en
The proof of 6.1. occupies the remaining point of this section. We always
assume the hypothesises of 6.1.

\en
\no
{\bf Lemma 6.2.} {\it
$W^{1,i}(E,g,h,\nabla) = W^{1,i}(E,g',h',\nabla')$, $0 \le i \le r$ as  
equivalent Banach spaces.                                    
}
                       
\en
\no
{\bf Corollary 6.3.} {\it
$W^{2,j}(E,g,h,\nabla) = W^{2,j}(E,g',h',\nabla')$, $0 \le j \le 
\frac{r}{2}$ as
equivalent Hilbert spaces. In particular,
\be
  L_2((M,E),g,h) = L_2((M,E),g',h').
\ee
}

\en
\no
{\bf Corollary 6.4.} 
$ H^j(E,D) \cong H^j(E,D'), \en 0 \le j \le \frac{r}{2}$.

\en
\no
{\bf Proof of (6.2).} This is well known for $h=h', \nabla' \in comp(\nabla)$.
But concerning $\nabla, \nabla'$ and $h, h'$ the only two facts needed in the
proof are just (6.3) (which is reformulated as (6.5)), (6.4) and the
equivalence of pointwise norms. The latter follows from (6.1), (6.2).
Into higher derivatives enter $(\nabla^g)^i, (\nabla^{g'})^j, i,j \le r-1$.
The conditions
\[ |\nabla^g-\nabla^{g'}|_{g,1,r-1} < \infty  \qu  
   |\nabla^g-\nabla^{g'}|_{g',1,r-1} < \infty     \]                             
follow from $g' \in comp(g)$. \hfill $\Box$

\en
\no
6.2. has a parallel version for the endomorphism bundle End $E$.

\en
\no
{\bf Lemma 6.5.} 
$\Omega^{1,1,i} (\mbox{End} \, E, g, h, \nabla) \cong 
\Omega^{1,1,i} (\mbox{End} \, E, g', h', \nabla'), \en 0 \le i \le r. $ 
\hfill $\Box$  

\en 
\no
{\bf Corollary 6.6.} 
$\Omega^{1,2,j} (\mbox{End} \, E, g, h, \nabla) \cong 
\Omega^{1,2,j} (\mbox{End} \, E, g', h', \nabla'), \en 0 \le j \le 
\frac{r}{2}$ .  
\hfill $\Box$  

\en
We obtain 
\be
  e^{-tD^2}, e^{-tD'^2}: L_2((M,E),g,h) \rightarrow H^j(E,D), 
  \en  0 \le j \le \frac{r}{2}
\ee
Hence 
\[ e^{-tD^2} - e^{-tD'^2} \]
is well defined. Our next task is to obtain an explicit expression for
$ e^{-tD^2} - e^{-tD'^2} $. For this we must modify Duhamel's principle
slightly. The steps 1. -- 5. in the proof of 5.2. remain unchanged.
We perform them for 
$(\cdot, \cdot)_q = h_q(\cdot, \cdot)$, 
$dq = dvol_q(g) \equiv dq(g)$. Then 5. reads as
\bea
   - \int\limits^\beta_\alpha \int\limits_M 
   h_q \Big( W(\tau,m,q), \big( D^2+\frac{\partial}{\partial t} \big)
   W'(t-\tau,q,p) \Big) dq(g) d\tau = \nonumber \\
   = \int\limits_M
   \Big[ h_q \big( W(\beta,m,q), W'(t-\beta,q,p) \big)- \nonumber \\
   - h_q \big( W(\alpha,m,q), W'(t-\alpha,q,p) \big) \Big] dq(g) . 
\eea  
Performing $\alpha \rightarrow 0^+, \beta \rightarrow t^-$ in (6.9)
yields
\bea  
   - \int\limits^t_0 \int\limits_M 
   h_q \Big( W(s,m,q), \big( D^2+\frac{\partial}{\partial t} \big)
   W'(t-s,q,p) \Big) dq(g) ds = \nonumber \\ 
   = \lim_{\beta \rightarrow t^-} \int\limits_M
   h_q \big( W(\beta,m,q), W'(t-\beta,q,p) \big) dq(g) 
   - W'(t,m,p) . 
\eea  
$g' \in comp(g)$ implies
\be
   dq(g) = \alpha(q) dq(g')
\ee
$0 < c_1 \le \alpha(q) \le c_2$. We rewrite
\bea
   &{}& \lim_{\beta \rightarrow t^-} \int\limits_M
   h_q \big( W(\beta,m,q), W'(t-\beta,q,p) \big) dq(g) \en = \nonumber \\
   &=& \lim_{\beta \rightarrow t^-} \int\limits_M  
   h' \big( W(\beta,m,q), W'(t-\beta,q,p) \big) \alpha(q)
   (\alpha(q)^{-1} dq(g)) \en + \nonumber \\   
   &+& \lim_{\beta \rightarrow t^-} \int\limits_M                       
   (h-h')_q \big( W(\beta,m,q), W'(t-\beta,q,p) \big) dq(g) 
   \en = \nonumber \\
   &=& \alpha(p) \cdot W(t,m,p) + 
   \lim_{\beta \rightarrow t^-} \int\limits_M
   (h-h')_q \big( W(\beta,m,q), W'(t-\beta,q,p) \big) dq(g) 
   \nonumber 
\eea
and obtain
\bea  
   &{}& - \int\limits^t_0 \int\limits_M   
   h_q \Big( W(s,m,q), \big( D^2+\frac{\partial}{\partial t} \big)
   W'(t-s,q,p) \Big) dq(g) ds \en = \nonumber \\   
   &=& - \int\limits^t_0 \int\limits_M 
   h_q \Big( W(s,m,q), \big( D^2-D'^2 \big)
   W'(t-s,q,p) \Big) dq(g) ds \en = \nonumber \\    
   &=& \alpha(p) W(t,m,p) - W'(t,m,p) \en + \nonumber \\
   &+& \lim_{\beta \rightarrow t^-} \int\limits_M  
   (h-h')_q \big( W(\beta,m,q), W'(t-\beta,p,q) \big) dq(g).
\eea
We see immediately that (6.12) expresses the operator equation
\bea
   e^{-tD^2}-e^{-tD'^2} = - \int\limits^t_0 e^{-sD^2}
   (D^2-D'^2) e^{-(t-s)D'^2} ds - \nonumber \\
   - \int\limits_M h'_p \Big( \lim_{\beta \rightarrow t^-}
   \int\limits_M (h-h')_q \big(W(\beta,m,q), W'(t-\beta,q,p) \big)
   dq(g), \cdot \Big) dp(g')
\eea
in $L_2((M,E),h',g')$ at kernel level. We want to show that both terms
on the right hand side of (6.13) are trace class operators with
uniformly bounded trace norm on compact $t$--intervalls 
$[a_0,a_1], a_0>0$, and we start with
\bea
   &{}& \int\limits^t_0 e^{-sD^2} (D^2-D'^2) e^{-(t-s)D'^2} ds \en =
   \nonumber \\                                                 
  &=& \int\limits^t_0 e^{-sD^2} D(D-D') e^{-(t-s)D'^2} ds \en + \\
  &+& \int\limits^t_0 e^{-sD^2} (D-D')D' e^{-(t-s)D'^2} ds .
\eea                                        
Write $D=\sum\limits^n_{i=1} e_i \cdot \nabla_{e_i} $ ,
$D'=\sum\limits^n_{i=1} e'_i \cdot' \nabla_{e_i'} $.
Then $(D-D') \Phi = \sum\limits^n_{i=1} e_i \cdot \nabla_{e_i} \Phi -
 e'_i \cdot' \nabla_{e_i'} \Phi $ . Consider 
$e \cdot \nabla_e - e' \cdot' \nabla_e' = (e-e') \cdot \nabla_e +
e' \cdot (\nabla_e - \nabla_{e'}) + e' \cdot (\nabla_{e'}- \nabla'_{e'})+
e'(\cdot - \cdot') \nabla'_{e'}$. Hence
\[ (D-D') \Phi = (\eta_1+\eta_2+\eta_3+\eta_4) \Phi , \]
where

$\eta_1 \Phi = \sum\limits_i (e_i-e'_i) \cdot \nabla_{e_i} \Phi$ ,
   
$\eta_2 \Phi = \sum\limits_i e'_i \cdot (\nabla_{e_i}- \nabla_{e'_i}) \Phi$ ,
                                                                    
$\eta_3 \Phi = \sum\limits_i e'_i \cdot (\nabla_{e'_i}- \nabla'_{e'_i}) \Phi$ ,
   
$\eta_4 \Phi = \sum\limits_i e'_i (\cdot -  \cdot') \nabla'_{e'_i} \Phi$ .
                       
\no                        
We simply write $\eta_\nu$ instead $\eta^{op}_\nu$ and obtain
\bea
   \mbox{(6.14) + (6.15)} &=& 
   \int\limits^t_0 e^{-sD^2} D (\eta_1+\eta_2+\eta_3+\eta_4)                       
   e^{-(t-s)D'^2} ds \en + \\
   &+& \int\limits^t_0 e^{-sD^2} (\eta_1+\eta_2+\eta_3+\eta_4)                       
   D' e^{-(t-s)D'^2} ds.  
\eea                           
We have to estimate
\be
  \int\limits^t_0 e^{-sD^2} D \eta_\nu e^{-(t-s)D'^2} ds 
\ee                                                     
and
\be
  \int\limits^t_0 e^{-sD^2} \eta_\nu D' e^{-(t-s)D'^2} ds .
\ee  
Decompose $\int\limits^t_0 = \int\limits^{\frac{t}{2}}_0 + 
\int\limits^t_{\frac{t}{2}}$ which yields
\def\theequation{$I_{\nu,1}$}   
\be
   \int\limits^{\frac{t}{2}}_0 e^{-sD^2} D \eta_\nu e^{-(t-s)D'^2} ds ,
\ee
\def\theequation{$I_{\nu,2}$}  
\be
   \int\limits^{\frac{t}{2}}_0 e^{-sD^2} \eta_\nu D' e^{-(t-s)D'^2} ds ,
\ee 
\def\theequation{$I_{\nu,3}$}  
\be
   \int\limits^t_{\frac{t}{2}} e^{-sD^2} D \eta_\nu  e^{-(t-s)D'^2} ds ,
\ee 
\def\theequation{$I_{\nu,4}$}  
\be
   \int\limits^t_{\frac{t}{2}} e^{-sD^2} \eta_\nu D' e^{-(t-s)D'^2} ds .
\ee                         

\def\theequation{\thesection.\arabic{equation}}  
\setcounter{equation}{19} 

\no                           
$(I_{\nu,j})$ has the same structure as $(I_j)$ in section 5.
                              
But in distinction to section 5, not all $\eta_\nu=\eta^{op}_\nu$                               
are operators of order zero. Only $\eta_3$ is zero order operator,
generated by an End $E$ valued 1--form $\eta_3$. For it we want to show
and then to use
\be
  |\eta_3|_{1,r-1} < \infty,
\ee
where $|\en|_{1,r-1}=|\en|_{g,h,\nabla,1,r-1}$ or $|\en|_{g',h',\nabla',1,r-1}$ as
we want. $\eta^{op}_1, \eta^{op}_2, \eta^{op}_4$ are first oder operators.
For them we want to show that their coefficients decrease sufficiently fast,
i. e. have finite $|\en|_{1,r-1}$--norm.

Altogether we have to estimate 16 integrands which split into even more.

We start with $\nu=3$. (6.20) is an immediate consequence of (6.3),(6.4)
and we are from an analytical point of view exactly in the situation of 
section 5. $(I_{3,1})-(I_{3,4})$ can be estimated quite parallel to 
$(I_1)-(I_4)$ in section 5 and we are done. There remains the estimate
of $(I_{\nu,j}), \nu \not= 3, j=1, \dots, 4$.
Start with $\nu=1, j=3 $ write
\bea
   &{}& e^{-sD^2} D \eta_1 e^{-(t-s)D'^2} = \nonumber \\
   &=& \Big( (D e^{-\frac{sD^2}{2}}) \circ f \Big) \circ 
   \Big( f^{-1} e^{-\frac{sD^2}{2}} \eta_1 \Big) \circ
   e^{-(t-s)D'^2} .
\eea
(6.21) holds since $e^{-\frac{sD^2}{2}}$  is a smoothing operator.
$e^{-(t-s)D'^2}$ is bounded in $[\frac{t}{2},t]$ and we handle and 
estimate it as in section 5. $ \Big( (D e^{-\frac{sD^2}{2}}) \circ f \Big)$ 
is Hilbert--Schmidt if $f \in L_2$. There remains to show that for
appropriate $f$
\be
   f^{-1} e^{-\frac{sD^2}{2}} \eta_1
\ee
is Hilbert--Schmidt. Recall $r>n+2, n \ge 2$, which implies
$\frac{r}{2} > \frac{n}{2} +1$,  
$r-1-n \ge \frac{r}{2} - \frac{n}{2}$,
$r-1 \ge \frac{r}{2}$, $2 \ge i$.
If we write in the sequel pointwise or Sobolev norms we should always
write $|\Psi|_{g',h',m}$, $|\Psi|_{H^r(E,D')}$, 
$|\Psi|_{g',h',\nabla',2,\frac{r}{2}}$, $|g-g'|_{g',m}$, 
$|g-g'|_{g',1,r}$ etc. But we omit the reference to
$g',h,\nabla',D,m'$ in the denotion for the sake of brevity. 
Moreover, as we already know, $g,h,\nabla,D$ generate equivalent norms.

Now $(\eta_1 \Phi)(m) = \sum\limits^n_{i=1} (e_i-e'_i) \cdot 
\nabla_{e_i} \Phi$,
\be
  |\eta_1 \Phi|_m \le \Big( \sum^n_{i=1} |e_i-e'_i|_m^2 
  \Big)^{\frac{1}{2}} \cdot |\nabla \Phi|_m \le 
  |g-g'|_m \cdot |\nabla \Phi|_m .
\ee 
Similarly, for supp $\Phi$ compact, $|\Phi|_{L_2}=1$,  
$|\Phi|_{H^{\frac{r}{2}}} \le C_1$,  $s>0$,
\bea
  &{}& |\eta_1 \Phi|_{H^{\frac{r}{2}-1}, B_{2d+s}(m)} \en \le \nonumber \\
  &\le& C_2 |g-g'|_{2,\frac{r}{2},B_{2d+s}(m)} \en \le 
  \qu ( \mbox{since} \en r>n+2 )  \nonumber \\
  &\le& C_3 |g-g'|_{1,r-1,B_{2d+s}(m)} \en  = \nonumber \\
  &=& C_3 \en vol(B_{2d+s}(m)) \cdot \nonumber \\ 
  &\cdot& \Big( \frac{1}{vol B_{2d+s}(m)}
  |g-g'|_{1,r-1,B_{2d+s}(m)} \Big) ,
\eea
and then we proceed as in (5.31)--(5.40), i. e. we set
\bea
  &{}& \tilde{\eta}_1(m) := C_4 \int\limits^\infty_0
  e^{-\frac{9}{10}\frac{s^2}{4t}}
  \Big( \frac{1}{vol B_{2d+s}(m)} 
  |g-g'|_{1,r-1,B_{2d+s}(m)} \Big) ds, \nonumber \\
  &{}& f_1(m) := \big( \tilde{\eta}_1(m) \big)^{\frac{1}{2}} , \nonumber             
\eea
and obtain finally
\bea
  &{}& |\mbox{integral kernel} \en f^{-1}_1 e^{-\frac{s}{2}D^2} 
  \eta_1 |_{L_2(M \times M)} \le C_5(t) , \nonumber \\
  &{}& \Big| \int\limits^t_{\frac{t}{2}} e^{-sD^2} D \eta_1  
  e^{-(t-s)D'^2} ds \Big|_1 \le C_6(\frac{t}{2},t) \cdot \sqrt{t} .
\eea
$(I_{1,3})$ is done. $(I_{1,1}), (I_{1,2}), (I_{1,4})$ can be handled
parallel to $(I_1), (I_2), (I_4)$ of section 5. If at the ''continuous''
end appear additional second derivatives, we proceed as with (5.54) using
the version of (5.37), i. e. (6.24).

Now it is completely clear that $(I_{\nu,j}), \nu=2,4, j= 1, \dots, 4$,
are done if we have an estimate for $\eta_2, \eta_4$ as above, coming
from our assumptions.
\bea                                                          
 |\eta_2 \Phi|_m &=& |\sum_i e'_i (\nabla_{e_i} - \nabla_{e'_i}) \Phi |_m       
 \en = \en |\sum_i e'_i \nabla_{e_i-e'_i} \Phi |_m \en \le \nonumber \\
 &\le& \Big( \sum_i |e_i-e'_i|^2_m \Big)^{\frac{1}{2}} |\nabla \Phi|_m 
 \en \le \en |g-g'|_m |\nabla \Phi|_m \nonumber 
\eea

\no        
Similarly for higher derivatives and we proceed as for $\eta_1$.
          
There remains $\eta_4$.
\bea                                                           
  (\eta_4 \Phi)(m) &=& \sum_i e'_i (\cdot - \cdot') \nabla'_{e'_i} \Phi 
  \en = \en \sum_i (\cdot - \cdot') (e'_i \otimes \nabla'_{e'_i} \Phi )
  \nonumber \\
  |\eta_4 \Phi|_m &\le& |\cdot - \cdot'|_m \cdot |\nabla' \Phi| \nonumber
\eea
Using our assumption (Clm)
\[ (\cdot - \cdot') \in W^{1,r} (\mbox{Hom} \en (TM \otimes E,E)), \]
we proceed as for the other $\eta_\nu$.

Finally we have to show that the operator
\be
  \Phi \rightarrow \int\limits_M h'_p \Big( \lim_{\beta \rightarrow t^-}
  \int\limits_M (h-h')_q \big( W(\beta,m,q), W'(t-\beta,q,p) \big) dq(g),
  \Phi(p) \Big) dp(g')
\ee
is a product of Hilbert-Schmidt operators. The first step is to rewrite
(6.26). For doing this, we apply the following facts

1. $\lim\limits_{\beta \rightarrow t^-} \int h'_p \big( 
W'(t-\beta,q,p), \cdot \big) dp(g') = \delta(q) $, 
$\Phi(p) \in C^\infty_c $

2. $W, W' \in C^\infty ({\bf R}_+ \times M \times M, 
E \Box \hspace{-0.42cm} \times E)$

3. For $a,b \in (V,(,)_V),a',b' \in (V',(,)_{V'}) $ holds 

$(a \otimes a',b \otimes b')_{V \otimes V'} = (a,b)_V \cdot (a',b')_{V'}
= ((a \otimes a',b),b') = (a,(a',b \otimes b')) $

4. The principle of Tonelli and the Fubini theorem for absolutely 
integrable integrands.

Then we can rewrite (6.26) as 
\be
  \Phi \rightarrow \int\limits_M (h-h')_q \big(W(t,m,q),\Phi(q) \big) dq(g)
\ee
We decompose (6.27) as
\bea
  &{}& \int\limits_M (h-h')_q \big( W(t,m,q), \Phi(q) \big) dq(g) 
  \en = \nonumber \\
  &=& \int\limits_M (h-h')_q \Big( \int\limits_M h_u 
  \big( W(\frac{t}{2},m,u), W(\frac{t}{2},m,u) \big) du(g),
  \Phi(q) \Big) dq(g) \en = \nonumber \\    
  &=& \int\limits_M h_u \Big(  W(\frac{t}{2},m,u) , 
  \int\limits_M (h-h')_q 
  \big( W(\frac{t}{2},u,q) , \Phi(q) \big) dq(g)
  \Big) du(g) \en = \nonumber \\    
  &=& A_2(A_1 \Phi), 
\eea 
where

$(A_1 \Phi)(u) = \int\limits_M (h-h')_q \big(  W(\frac{t}{2},u,q) , 
\Phi(q) \big) dq(g) $   

\no
and

$(A_2 \Psi)(m) = \int\limits_M h_u \big(  W(\frac{t}{2},m,u) , 
\Psi(u) \big) du(g) $ .

\no
Next we want to write 
\be
  A_2 \circ A_1 = (A_2 \circ f \cdot) \circ ((f^{-1} \cdot) \circ A_1 ), 
\ee
$f$ a scalar function s. t. $A_2 \circ f$ and $(f^{-1} \cdot) \circ A_1$
are Hilbert-Schmidt operators and we start with $A_1$. Our procedure is
as above. We estimate the integral norm of $A_1$ with respect to one
variable and then we define $f$. Now
\bea
  &{}& \Big| (h-h')_q \big( W(t,u,q), \cdot \big) 
  \Big|^2_{L_2(h,dq(g))} \en \le \nonumber \\
  &\le& \Big| h \big( W(t,u,q), (|h-h'|_{h,q} \cdot ) \cdot 
  \big) \Big|_{L_2(h,dq(g))} .
\eea
This amounts as in section 5 to estimate
\[ \sup_{
   \begin{array}{c} 
   \Phi \in C^\infty_c \\
   |\Phi|_{L_2} = 1 \\
   |\Phi|_{H^{\frac{r}{2}}} \le C
   \end{array}
   }        
   |N(\Phi)|^2 , \]
where
\bea
  N(\Phi) &=& \Big| \big\langle \delta(u), e^{-tD^2}
  (|h-h'|_{h,q} \cdot \Phi) \big\rangle \Big| \en = \nonumber \\
  &=& \frac{1}{\sqrt{4 \pi t}}
  \Big| \int\limits^{+ \infty}_{- \infty} e^{-\frac{s^2}{4t}}  
  e^{isD} (|h-h'|_h \cdot \Phi ) ds \Big| .
\eea
But now we proceed literally as is (5.19) -- (5.35), replacing $\eta$ by
$|h-h'|$ and setting
\[ \tilde{f} := C_{10} \cdot \int\limits^\infty_0 
   e^{-\frac{9}{10} \frac{s^2}{4} \frac{t}{2}} 
   \Big( \frac{1}{vol B_{2d+s}(u)} |h-h'|_{h,\nabla,r-1, B_{2d+s}(u)} 
   \Big) ds . \]
Then again (for $h \not= h'$) $\tilde{f}(u)$ is $\not= 0$ everywhere
and 1--summable,
\be
  \int\limits_M \Big| (h-h')_q  \big( W(\frac{t}{2},u,q), \cdot  
  \big) \Big|^2_{h,q} dq(g) \le \tilde{f}(u)^2 .
\ee
We set $f(u) := (\tilde{f}(u))^{\frac{1}{2}} $ and infer $f(u) \not= 0$
everywhere, $f \in L_2$ and
\bea
  &{}& \int\limits_M \int\limits_M f(u)^{-2} \Big|  
  (h-h')_q \big( W(\frac{t}{2},u,q), \cdot \big)
  \Big|^2 dq(g) du(g) \en \le \nonumber \\  
  &{}& \le \en \int\limits_M \frac{1}{\tilde{f}(u)} \tilde{f}(u)^2  du(g) 
  \en = \en \int\limits_M \tilde{f}(u) du(g) \en \le \en C_{11}(t) 
  \cdot \sqrt{t} ,
\eea
where for $t \in [a_0,a_1], a_0>0$,
\be
  \sup_{t \in [a_0,a_1]} C_{11} (t) \cdot \sqrt{t} = C_{12} (a_0,a_1) < \infty .
\ee
Quite similarly as in (5.12) -- (5.14), 
\be
  \int\limits_M \int\limits_M f(m)^2 \big| W(\frac{t}{2},m,u) 
  \big|^2_h du(g) dm(g) \le C_{13}(t) , 
\ee
and finally, as in (5.39) -- (5.40)
\bea
   & \Big| \Phi \rightarrow \int\limits_M h'_p \Big( 
   \lim_{\beta \rightarrow t^-} \int\limits_M (h-h')_q
   \big( W(\beta,m,q), W'(t-\beta,q,p) \big) dq(g),
   \Phi(p)                           
   \Big) dp(g') \Big|_1  \\ \nonumber &   
   \le C_{14}(t),                      
\eea                                   
where for $[a_0, a_1], a_0>0$,
\be
  \sup_{t \in [a_0, a_1]} C_{14}(t) = C_{15}(a_0,a_1) < \infty ,
\ee
this finishes the proof of 6.1. \hfill $\Box$

\en
\no
{\bf Example.} The simplest standard example is 
$E=(\Lambda^* T^* M, g_{\Lambda^*}, \nabla^{g_{\Lambda^*}})$ with Clifford
multiplication
\[ x \otimes \omega \in T_m M \otimes \Lambda^* T^*_m M \rightarrow 
   X \cdot \omega := \omega_X \wedge \omega - i_X \omega ,\]
where $ \omega_X := g(\cdot, X)$. In this case, $E$ as a vector bundle
remains fixed but the Clifford module structure varies smoothly with
$g, g' \in comp(g)$, i.e. (6.1), automatically implies (Clm), (6.2),
(6.3), (6.4). It is well known that in this case $D=d+d^*$, 
$D^2=(d+d^*)^2$ = graded Laplace operator $\bu$. Hence we obtain

\en
\no
{\bf Corollary 6.7.} {\it 
Assume $(M^n,g)$ with $(I), (B_k), k \ge r >n+2,$
$g' \in {\cal M}(I,B_k)$, $ g' \in comp(g) \subset {\cal M}^{1,r} (I, B_k)$.
Then 
\[ e^{-t \bu} - e^{-t \bu'} \]
is of trace class and the trace norm is uniformly bounded on compact
$t$--intervalls $[a_0,a_1], a_0>0$.} \hfill $\Box$

\en
\no
{\bf Remark.} 
We are also able to prove 6.7 directly without reference to 6.1. For
this we write $\bu'=\bu+\eta$, calculate and estimate $\eta$ (which
is very easy), apply Duhamel's principle and proceed as before.
\hfill $\Box$                                                 

\me
We need in section 7 the theorem analoguos to 5.3 for the case of 
additional variation of the metrics.  

\me
\no
{\bf Theorem 6.8.} {\it Suppose the hypothesises of 6.1, replacing
$r>n+2$ by $r>n+3$. Then 
\[ De^{-tD^2} - D'e^{-tD'^2} \]
is of trace class and the trace norm ist uniformly bounded on 
compact $t$--intervalls $[a_0,a_1], a_0>0$. }

\me
\no
{\bf Proof.} The proof is a simple combination of the proofs of
5.3 and 6.1. \hfill $\Box$

\section{Relative index theory}
\setcounter{equation}{0}

We now assume that $E$ is endowed with an involution 
$\tau : E \rightarrow E$, s. t.
\bea
   &{}& \tau^2 = 1, \tau^* = \tau ,\\
   &{}& [\tau,X]_+ = 0 \en \mbox{for} \en X\in TM , \\
   &{}& [\nabla,\tau] = 0 .
\eea
Then $L_2(M,E) = L_2(M,E^+) \oplus L_2(M,E^-)$,
\[ D = 
   \left( \begin{array}{cc} 
   0 & D^- \\
   D^+ & 0  
   \end{array} \right) \]
and $D^- = (D^+)^*$. If $M^n$ is compact then as usual
\be
  ind \, D := ind \, D^+ := dim \, ker \, D^+ - dim \, ker \, D^-
  \equiv tr \, (\tau e^{-tD^2}) ,
\ee
where we unterstand $\tau$ as 
\[ \tau =
   \left( \begin{array}{cc} 
   I & 0 \\
   0 & -I  
   \end{array} \right) :
   L_2(E^*) \oplus L_2(E^-) \rightarrow L_2(E^+) \oplus L_2(E^-) . \]
For open $M^n$ $ind \, D$ in general is not defined since
$dim \, ker \, D^+$, $dim \, ker \, D^-$ are not of trace class.

The appropriate approach on open manifolds is relative index theory
for pairs of operators $D, D'$. If $e^{-tD^2} - e^{-tD'^2}$ is of
trace class then
\be
  ind \, (D, D') := tr \, ( \tau ( e^{-tD^2} - e^{-tD'^2} ))   
\ee
makes sense, but at the first glance (7.5) should depend on $t$.

\me
\no                    
{\bf Proposition 7.1.} {\it
Suppose $e^{-tD^2} - e^{-tD'^2}$ and $D e^{-tD^2} - D' e^{-tD'^2}$
of trace class for all $t>0$ and $|D e^{-tD^2} - D' e^{-tD'^2}|_1$  
uniformly bounded on any compact t--intervall $[a_0,a_1], a_0>0$.
Then $tr \, ( \tau ( e^{-tD^2} - e^{-tD'^2} ))$ is independent of $t$.
}
                       
\me
\no
See [1] for a proof. \hfill $\Box$

\me
\no                    
{\bf Corollary 7.2.} {\it 
Assume the hypotheses of 7.1. Then $ind \, (D,D')$ is independent of t and
hence well defined.}         

\me
\no                    
{\bf Corollary 7.3.} {\it 
Assume the hypotheses of 7.1 and $inf \, \sigma_e (D^2) >0 $. Then 
$ind \, D, \, ind \, D'$ are well defined and 
\be
  tr \, ( \tau ( e^{-tD^2} - e^{-tD'^2} )) = ind \, D - ind \, D' .
\ee
} 

\me
\no
{\bf Proof.} 
From our assumptions, $\sigma_e(D^2)=\sigma_e(D'^2)$.
$\inf \sigma_e(D^2)>0$ immediately implies
$\dim \, ker \, D^+, \, dim \, ker \, D^- < \infty $.
According to (0.7) of [11], there exists a constant $c>0$ s. t.
\be
  tr \, ( \tau ( e^{-tD^2} - e^{-tD'^2} )) = ind \, D - ind \, D'
  + O(e^{-ct}) .
\ee
Performing $\lim\limits_{t \rightarrow \infty}$ in (7.7) and using 7.2,
we obtain (7.6). \hfill $\Box$

\me
Assume now the hypotheses of 5.1. Then we have asymptotic expansions
\be
  tr \, \big( \tau W(t,m,m) \big)
  \begin{array}{c}
   {} \\ \sim \\ t \rightarrow 0^+
  \end{array}
  t^{-\frac{n}{2}} b_{-\frac{n}{2}} (D,m) + \dots + b_0 (D,m) +
  b_{\frac{1}{2}} (D,m) t^{\frac{1}{2}} + \dots
\ee  
\be
  tr \, \big( \tau W'(t,m,m) \big)
  \begin{array}{c}
   {} \\ \sim \\ t \rightarrow 0^+
  \end{array}   
  t^{-\frac{n}{2}} b_{-\frac{n}{2}} (D',m) + \dots + b_0 (D',m) +
  b_{\frac{1}{2}} (D',m) t^{\frac{1}{2}} + \dots
\ee  
We show in the next section that 
\be
  b_i (D,m) - b_i (D',m) \in L_1 , \qu 
  - \frac{n}{2} \le i \le 1 .
\ee
Define
\[ ind_{top} \, (D,D') := \int\limits_M 
   \big( b_0(D,m) - b_0 (D',m) \big) dm . \]
According to (7.10), $ind_{top} \, (D,D')$ is well defined.

\me
\no
{\bf Theorem 7.4.} {\it
Suppose the hypotheseses of 5.3. Then 
\[ ind \, (D,D') = ind_{top} \, (D,D') \]
If additionally $\inf \sigma_e (D^2) > 0$ then
\[ ind_{top} \, (D,D') = ind \, D - ind \, D' . \]
}

\me 
\no
{\bf Proof.}
This follows immediately from 7.2, 7.3, (7.8), (7.9) and the fact that
the $L_2$--trace of a trace class integral operator is given by the
integral of the kernel on the diagonal (after forming pointwise
traces). \hfill $\Box$

If we admit variation of $g$ too as in section 6, then the heat kernel of 
$e^{-tD'^2}$ in $L_2((M,E),g,h)$ is given by 
$\alpha(m)^{\frac{1}{2}} W'(t,m,p) \alpha(p)^{-\frac{1}{2}}$. But on the diagonal the
$\alpha$'s cancel out and the asymptotic expansion of $W'(t,m,m)$ with 
respect to $L_2(g)$ is the same as with respect to $L_2(g')$. We
obtain for $W(t,m,m)$ or $W'(t,m,m)$ heat kernel coefficients 
$b_i(D(g,h,\nabla),m)$ or $b_i(D(g',h',\nabla'),m)$, respectively.
We show in the next section that under the hypotheses of 6.1
\be 
  b_i(D(g,h,\nabla),m) - b_i(D(g',h',\nabla'),m) \in L_1 , 
  \qu -\frac{n}{2} \le i \le 1 .
\ee

\me
\no
{\bf Theorem 7.5} {\it
Suppose the hypotheseses of 6.8. Then
\[ ind \, (D,D') = ind_{top} \, (D,D') . \]
If additionally $\inf \sigma_e(D^2) > 0$, then
\[ ind_{top} \, (D,D') = ind \, D - ind \, D' . \]
}

\no
The proof runs through literally as that of 7.4. \hfill $\Box$

\me 
Deeper results on the relative index using scattering theory will be
established in a forthcoming paper.

\section{Relative $\zeta$--functions, determinants and torsion}
\setcounter{equation}{0}

We start with a pair $D, D'$ assuming the hypotheses of 5.1. Then we
have the asymptotic expansion
\be
  tr \, W(t,m,m) 
  \begin{array}{c}
   {} \\ \sim \\ t \rightarrow 0^+
  \end{array}    
  t^{-\frac{n}{2}} b_{-\frac{n}{2}} (m) + 
  t^{-\frac{n}{2}+1} b_{-\frac{n}{2}+1} + \dots
\ee
and analogously for $tr \, W'(t,m,m)$ with 
\[ b_{-\frac{n}{2}+1}(m) =  b_{-\frac{n}{2}+l} (D(g,h,\nabla),m), \qu 
   b'_{-\frac{n}{2}+l} (m) =  b_{-\frac{n}{2}+l} (D(g,h,\nabla'),m) . \]
The heat kernel coefficients have for $l \ge 1$ a representation
\be
  b_{-\frac{n}{2}+l} = \sum^l_{k=1} \sum^l_{q=0} 
  \sum_{i_1, \dots , i_k \ge 0} \nabla^{i_1} R^g \dots \nabla^{i_q} R^g
  tr \, (\nabla^{i_{q+1}} R^E \dots \nabla^{i_k} R^E ) C^{i_1, \dots , i_k} ,
\ee                                          
where $C^{i_1, \dots , i_k}$ stands for a contraction with respect to $g$, 
i.e. it is built up by linear combination of products of the $g^{ij}$.

\me
\no
{\bf Lemma 8.1.} {\it
$b_{-\frac{n}{2}+l} - b'_{-\frac{n}{2}+i} \in L_1(M,g)$, 
$0 \le l \le \frac{n+3}{2}$. }          

\me
\no
{\bf Proof.} Forming the difference 
$b_{-\frac{n}{2}+l} - b'_{-\frac{n}{2}+l}$, we obtain a sum of terms of
the kind 
\be
   \nabla^{i_1} R^g \dots \nabla^{i_q} R^g \, tr \, 
   [\nabla^{i_{q+1}} R^E \dots \nabla^{i_k} R^E - 
   \nabla'^{i_{q+1}} R'^E \dots \nabla'^{i_k} R'^E] .
\ee
$g$ is here fixed. The highest derivative of $R^q$ with respect to $\nabla^g$
occurs if $q=k, i_1= \dots = i_{q-1} = 0$. Then we have
\be
   (\nabla^g)^{2l-2k} . 
\ee
By assumption, we have bounded geometry of order $\ge r > n+2$,
i. e. of order $\ge n+3$. Hence $(\nabla^g)^i R^g$ is bounded for
$i \le n+1$. To obtain bounded $\nabla^j R^g$--coefficients of [ \dots ] in 
(8.3), we must assume
\be
  2l-2 \le n+1, \qu l \le \frac{n+3}{2} .
\ee
Similarly we see that the highest occuring derivatives of $R^E$, 
$R'^E$ in [ \dots ] are of order $2l-2$. The corresponding expression
\be
   R^E \nabla^{2l-2} R^E - R'^E \nabla'^{2l-2} R'^E =
   (R^E-R'^E)(\nabla^{2l-2} R^E) + R'^E 
   (\nabla^{2l-2} R^E - \nabla^{2l-2} R'^E) .
\ee
We want to apply the module structure theorem.
$\nabla - \nabla' \in \Omega^{1,1,r} ({\cal G}^{Cl}_E, \nabla) =
 \Omega^{1,1,r} ({\cal G}^{Cl}_E, \nabla')$ 
implies $R^E-R'^E \in \Omega^{2,1,r-1}$.
We can apply the module structure theorem (and conclude that all norm
products of derivatives of order $\le 2l-2$ are absolutely
integrable) if $2l-2 \le r-1$, $2l-2 \le n+1$, $l \le \frac{n+3}{2}$.
Hence, (8.5) $\in L_1$ since $R^E, R'^E$ bounded. It is now a very
simple combinatorial matter to write [ \dots ] in (8.3) as a sum of terms
each of them is a product of differences $(\nabla^i R^E - \nabla'^i R'^E)$
with bounded functions $\nabla^j R^E$, $\nabla'^{j'} R'^E$. Remember
$\nabla, \nabla' \in {\cal C}_E(B_k)$. This proves 8.1. \hfill $\Box$

\me
\no
{\bf Lemma 8.2.} {\it 
There is an expansion
\be
  tr \, (e^{-tD^2}-e^{-tD'^2})= t^{-\frac{n}{2}} a_{-\frac{n}{2}}+
  \dots + t^{-\frac{n}{2}+[\frac{n+3}{2}]} a_{-\frac{n}{2}+[\frac{n+3}{2}]}
  + O(t^{-\frac{n}{2}+[\frac{n+3}{2}]+1}) .
\ee
}

\me
\no
{\bf Proof.} Set
\be
  a_{-\frac{n}{2}+i} = \int \Big( b_{-\frac{n}{2}+i}(m) -
  b'_{-\frac{n}{2}+i}(m)\Big) dm
\ee
and use
\bea
  tr \, W(t,m,m) &=&  t^{-\frac{n}{2}} b_{-\frac{n}{2}}+
  \dots + t^{-\frac{n}{2}+[\frac{n+3}{2}]} b_{-\frac{n}{2}+[\frac{n+3}{2}]}
  + O(m, t^{-\frac{n}{2}+[\frac{n+3}{2}]+1}) , \\           
   tr \, W'(t,m,m) &=&  t^{-\frac{n}{2}} b'_{-\frac{n}{2}} + \dots 
  + O'(m, t^{-\frac{n}{2}+[\frac{n+3}{2}]+1}) \\
  tr \, (e^{-tD^2}-e^{-tD'^2}) &=& \int \Big( 
  tr \, W(t,m,m) - tr \, W'(t,m,m) \Big) dm . \nonumber
\eea                 
The only critical point is 
\be
  \int\limits_M O(m, t^{-\frac{n}{2}+[\frac{n+3}{2}]+1}) - 
  O'(m, t^{-\frac{n}{2}+[\frac{n+3}{2}]+1}) dm = 
  O(t^{-\frac{n}{2}+[\frac{n+3}{2}]+1}) .  
\ee                      
(8.11) requires a very careful investigation of the concrete
representatives for $O(m, t^{-\frac{n}{2}+[\frac{n+3}{2}]})$.
We did this step by step, following [9], p. 21/22, 50--51. Very
roughly spoken, the $m$--dependence of $O(m,\cdot)$ is given by 
the parametrix construction, i. e. by differences
of corresponding derivatives of the $\Gamma^\beta_{i \alpha}$,
$\Gamma'^\beta_{i \alpha}$, which are integrable by assumption.
\hfill $\Box$

\me
\no
{\bf Definition.} 
Assume the hypotheses of 5.1. Set
\be
  \zeta_1 (s,D,D') := \frac{1}{\Gamma(s)} \int\limits^1_0 t^{s-1}
  tr \, (e^{-tD^2}-e^{-tD'^2}) dt.
\ee                     
Using 8.7, 
\be
  \int\limits^1_0 t^{s-1} t^{-\frac{n}{2}+[\frac{n+3}{2}]} dt =
  \frac{1}{s-\frac{n}{2}+[\frac{n+3}{2}]} ,
\ee
\be
  \frac{1}{\Gamma(s)} \int\limits^1_0 t^{s-1} 
  O(t^{-\frac{n}{2}+[\frac{n+3}{2}]+1}) dt \en \mbox{holomorphic for} \en
  Re(s)+(-\frac{n}{2})+[\frac{n+3}{2}]+1>0
\ee
and $[\frac{n+3}{2}] \ge \frac{n}{2}+1$, we
obtain a function meromorphic in $Re(s) > -1$, holomorphic in 
$s=0$ with simple poles at $s=\frac{n}{2}-l$, $l \le [\frac{n+3}{2}]$.
Assume additionally $\inf \sigma_e (D^2) > 0$ and set
$h= dim \, ker \, D^2 - dim \, ker \, D'^2 $. Then 
\be
 tr \, (e^{-tD^2}-e^{-tD'^2}) =   dim \, ker \, D^2 - dim \, ker \, D'^2
 + O(e^{-ct}) = h + O(e^{-ct}) \en \mbox{for} \en 
 t \rightarrow \infty, \en c>0 .
\ee
Define for $Re(s) < 0$
\[ \zeta_2(s,D^2,D'^2) := \frac{1}{\Gamma(s)} \int\limits^\infty_1 
   t^{s-1} tr \, (e^{-tD^2}-e^{-tD'^2}) dt =    
   \frac{1}{\Gamma(s)} \int\limits^\infty_1 t^{s-1}(h+O(e^{-ct})) dt .  \]
$\zeta_2(s,D^2,D'^2)$ is holomorphic in $Re(s)<0$ and admits a 
meromorphic extension to ${\bf C}$ which is holomorphic in $s=0$.

\no
Define
\[ \zeta(s,D^2,D'^2) := \zeta_1(s,D^2,D'^2) + \zeta_2(s,D^2,D'^2) .\]
We proved the following

\me        
\no
{\bf Theorem 8.3.} {\it
Suppose the hypotheseses of 5.1 and additionally $\inf \sigma _e(D^2)>0$.
Then 
\[ \zeta(s,D^2,D'^2) = \frac{1}{\Gamma(s)} \int\limits^\infty_0  
   t^{s-1} tr \, (e^{-tD^2}-e^{-tD'^2}) ds   \]  
is after meromorphic extension well defined in $Re(s)>-1$ and holomorphic    
in $s=0$.} \hfill $\Box$

\me
\no
{\bf Definition.} Suppose the hypotheses of 8.3. Then
\[ \mbox{det} \, (D^2,D'^2) := e^{-\zeta'(0,D^2,D'^2)} \]  
is well defined and is called the relative determinant of $D^2 / D'^2$.

\me
\no
{\bf Remark.} Fix $g,h,\nabla_0,D_0 = D(g,h,\nabla_0)$. Then we defined for any
$D(g,h,\nabla)$, $\nabla \in comp(\nabla_0) \cap {\cal C}_E(B_k) \subset 
{\cal C}^{1,r}_E(B_k) $ the relative determinant $\mbox{det} \, (D^2,D'^2)$.
\hfill $\Box$

If we suppose the hypotheses of 6.1 then we can repeat the preceding
considerations and estimates word by word. 

\me
\no
{\bf Theorem 8.4.} {\it
Suppose the hypotheses of 6.1. and $\inf \sigma_e(D^2)>0$.
Then 
\[ \zeta(s,D^2,D'^2) = \frac{1}{\Gamma(s)} \int\limits^\infty_0  
   t^{s-1} tr \, (e^{-tD^2}-e^{-tD'^2}) dt    \]  
is after meromorphic extension well defined in $Re(s)>-1$ and 
holomorphic in $s=0$. Hence the relative determinant
\[ \mbox{det} \, (D^2,D'^2) := e^{-\zeta'(0,D^2,D'^2)} \]  
is well defined. } \hfill $\Box$

\me
\no
{\bf Corollary 8.5.} {\it
Suppose $(M^n,g)$ with $(I)$, $(B_k)$, $k \ge r > n+2$, 
$g' \in comp(g) \cap {\cal M}(I,B_k) \subset {\cal M}^{1,r}(I,B_k)$,
and additionally $\inf \sigma_e (\bu_q) > 0$, $q=0, \dots, n$.
Then the relative analytic torsion $\tau(M,g,g')$,
\[ \log \tau(M,g,g') := \frac{1}{2} \sum^n_{q=0} (-1)^q 
   \frac{d}{ds} \zeta (s,\bu_q,\bu'_q) \Big|_{s=0} \]
is well defined.} \hfill $\Box$

\me
In a forthcoming paper we drop considerably the assumption
$\inf \sigma_e (\cdot)>0$ and discuss further applications.

\end{document}